\documentclass[10pt]{amsart}
\usepackage{amsfonts}
\usepackage{mathrsfs}
\usepackage{amscd}
\usepackage{amsmath}
\usepackage{amssymb}
\usepackage{latexsym}
\usepackage{lscape}
\usepackage{xypic}
\usepackage{comment}
\usepackage{amscd}
\usepackage{wasysym}  
\usepackage{tikz}
\usetikzlibrary{matrix,arrows}
\usepackage{tikz-cd}
\usepackage{appendix}
\usepackage{geometry}
\usepackage[nice]{nicefrac}
\usepackage{dsfont, stmaryrd}

\usepackage[pagebackref,hyperindex,linktocpage=true]{hyperref}
\hypersetup{
    colorlinks,
    linkcolor={red!50!black},
    citecolor={blue!50!black},
    urlcolor={blue!80!black}
}

\usepackage{color}

\usetikzlibrary{decorations.markings}

\makeatletter
\tikzcdset{
  open/.code     = {\tikzcdset{hook, circled};},
  closed/.code   = {\tikzcdset{hook, slashed};},
  open'/.code    = {\tikzcdset{hook', circled};},
  closed'/.code  = {\tikzcdset{hook', slashed};},
  circled/.code  = {\tikzcdset{markwith = {\draw (0,0) circle (.375ex);}};},
  slashed/.code  = {\tikzcdset{markwith = {\draw[-] (-.4ex,-.4ex) -- (.4ex,.4ex);}};},
  markwith/.code ={
    \pgfutil@ifundefined%
    {tikz@library@decorations.markings@loaded}%
    {\pgfutil@packageerror{tikz-cd}{You need to say %
      \string\usetikzlibrary{decorations.markings} to use arrows with markings}{}}{}%
    \pgfkeysalso{/tikz/postaction = {
      /tikz/decorate,
      /tikz/decoration={markings, mark = at position 0.5 with {#1}}}
    }
  },
}
\makeatother

\renewcommand{\geq}{\geqslant}
\renewcommand{\leq}{\leqslant}
\renewcommand{\ge}{\geqslant}
\renewcommand{\le}{\leqslant}

\newtheorem{thm}{Theorem}[section]

\newtheorem{propo}[thm]{Proposition}
\newtheorem{lem}[thm]{Lemma}
\newtheorem{sublem}[thm]{Sublemma}
\newtheorem{lem-def}[thm]{Lemma-Definition}
\newtheorem{cor}[thm]{Corollary}
\newtheorem{conject}[thm]{Conjecture}
\newtheorem{propert}[thm]{Properties}
\newtheorem{observ}[thm]{Observation}
\newtheorem{assum}[thm]{Assumption}

\theoremstyle{definition}
\newtheorem*{ack}{Acknowledgement}

\newtheorem{ex}[thm]{Example}

\newtheorem{rmk}[thm]{Remark}
\newtheorem{dfn}[thm]{Definition}
\newtheorem{quest}[thm]{Question}

\numberwithin{equation}{section}

\newcommand{\nc}{\newcommand}

\nc{\theo}{\begin{thm}} \nc{\xtheo}{\end{thm}}
\nc{\prop}{\begin{propo}} \nc{\xprop}{\end{propo}}
\nc{\lemm}{\begin{lem}} \nc{\xlemm}{\end{lem}}
\nc{\sublemm}{\begin{sublem}} \nc{\xsublemm}{\end{sublem}}
\nc{\lemmdefi}{\begin{lem-def}} \nc{\xlemmdefi}{\end{lem-def}}
\nc{\coro}{\begin{cor}} \nc{\xcoro}{\end{cor}}
\nc{\conj}{\begin{conject}} \nc{\xconj}{\end{conject}}
\nc{\proper}{\begin{propert}} \nc{\xproper}{\end{propert}}
\nc{\obse}{\begin{observ}} \nc{\xobse}{\end{observ}}
\nc{\ques}{\begin{quest}} \nc{\xques}{\end{quest}}

\nc{\ackn}{\begin{ack}} \nc{\xackn}{\end{ack}}
\nc{\exam}{\begin{ex}} \nc{\xexam}{\end{ex}}
\nc{\rema}{\begin{rmk}} \nc{\xrema}{\end{rmk}}
\nc{\defi}{\begin{dfn}} \nc{\xdefi}{\end{dfn}}

\nc{\pf}{\begin{proof}} \nc{\xpf}{\end{proof}}

\nc{\on}{\operatorname}
\nc{\fraka}{{\mathfrak a}} \nc{\bba}{{\mathbf a}}
\nc{\frakb}{{\mathfrak b}}
\nc{\frakc}{{\mathfrak c}}
\nc{\frakd}{{\mathfrak d}}
\nc{\frake}{{\mathfrak e}}
\nc{\frakf}{{\mathfrak f}}
\nc{\frakg}{{\mathfrak g}}
\nc{\frakh}{{\mathfrak h}}
\nc{\fraki}{{\mathfrak i}}
\nc{\frakj}{{\mathfrak j}}
\nc{\frakk}{{\mathfrak k}}
\nc{\frakl}{{\mathfrak l}}
\nc{\frakm}{{\mathfrak m}}
\nc{\frakn}{{\mathfrak n}}
\nc{\frako}{{\mathfrak o}}
\nc{\frakp}{{\mathfrak p}}
\nc{\frakq}{{\mathfrak q}}
\nc{\frakr}{{\mathfrak r}}
\nc{\fraks}{{\mathfrak s}}
\nc{\frakt}{{\mathfrak t}}
\nc{\fraku}{{\mathfrak u}}
\nc{\frakv}{{\mathfrak v}}
\nc{\frakw}{{\mathfrak w}}
\nc{\frakx}{{\mathfrak x}}
\nc{\fraky}{{\mathfrak y}}
\nc{\frakz}{{\mathfrak z}}
\nc{\frakA}{{\mathfrak A}}
\nc{\frakB}{{\mathfrak B}}
\nc{\frakC}{{\mathfrak C}}
\nc{\frakD}{{\mathfrak D}}
\nc{\frakE}{{\mathfrak E}}
\nc{\frakF}{{\mathfrak F}}
\nc{\frakG}{{\mathfrak G}}
\nc{\frakH}{{\mathfrak H}}
\nc{\frakI}{{\mathfrak I}}
\nc{\frakJ}{{\mathfrak J}}
\nc{\frakK}{{\mathfrak K}}
\nc{\frakL}{{\mathfrak L}}
\nc{\frakM}{{\mathfrak M}}
\nc{\frakN}{{\mathfrak N}}
\nc{\frakO}{{\mathfrak O}}
\nc{\frakP}{{\mathfrak P}}
\nc{\frakQ}{{\mathfrak Q}}
\nc{\frakR}{{\mathfrak R}}
\nc{\frakS}{{\mathfrak S}}
\nc{\frakT}{{\mathfrak T}}
\nc{\frakU}{{\mathfrak U}}
\nc{\frakV}{{\mathfrak V}}
\nc{\frakW}{{\mathfrak W}}
\nc{\frakX}{{\mathfrak X}}
\nc{\frakY}{{\mathfrak Y}}
\nc{\frakZ}{{\mathfrak Z}}
\nc{\bbA}{{\mathbb A}}
\nc{\bbB}{{\mathbb B}}
\nc{\bbC}{{\mathbb C}}
\nc{\bbD}{{\mathbb D}}
\nc{\bbE}{{\mathbb E}}
\nc{\bbF}{{\mathbb F}} \nc{\bbf}{{\mathbf f}}
\nc{\bbG}{{\mathbb G}}
\nc{\bbH}{{\mathbb H}}
\nc{\bbI}{{\mathbb I}}
\nc{\bbJ}{{\mathbb J}}
\nc{\bbK}{{\mathbb K}}
\nc{\bbL}{{\mathbb L}}
\nc{\bbM}{{\mathbb M}}
\nc{\bbN}{{\mathbb N}}
\nc{\bbO}{{\mathbb O}}
\nc{\bbP}{{\mathbb P}}
\nc{\bbQ}{{\mathbb Q}}
\nc{\bbR}{{\mathbb R}}
\nc{\bbS}{{\mathbb S}}
\nc{\bbT}{{\mathbb T}}
\nc{\bbU}{{\mathbb U}}
\nc{\bbV}{{\mathbb V}}
\nc{\bbW}{{\mathbb W}}
\nc{\bbX}{{\mathbb X}}
\nc{\bbY}{{\mathbb Y}}
\nc{\bbZ}{{\mathbb Z}}
\nc{\calA}{{\mathcal A}}
\nc{\calB}{{\mathcal B}}
\nc{\calC}{{\mathcal C}}
\nc{\calD}{{\mathcal D}}
\nc{\calE}{{\mathcal E}}
\nc{\calF}{{\mathcal F}}
\nc{\calG}{{\mathcal G}}
\nc{\calH}{{\mathcal H}}
\nc{\calI}{{\mathcal I}}
\nc{\calJ}{{\mathcal J}}
\nc{\calK}{{\mathcal K}}
\nc{\calL}{{\mathcal L}}
\nc{\calM}{{\mathcal M}}
\nc{\calN}{{\mathcal N}}
\nc{\calO}{{\mathcal O}}
\nc{\calP}{{\mathcal P}}
\nc{\calQ}{{\mathcal Q}}
\nc{\calR}{{\mathcal R}}
\nc{\calS}{{\mathcal S}}
\nc{\calT}{{\mathcal T}}
\nc{\calU}{{\mathcal U}}
\nc{\calV}{{\mathcal V}}
\nc{\calW}{{\mathcal W}}
\nc{\calX}{{\mathcal X}}
\nc{\calY}{{\mathcal Y}}
\nc{\calZ}{{\mathcal Z}}

\nc{\scrA}{{\mathscr A}}
\nc{\scrE}{{\mathscr E}}
\nc{\scrR}{{\mathscr R}}

\nc{\Bmu}{\mbox{$\raisebox{-0.59ex}{$l$}\hspace{-0.18em}\mu\hspace{-0.88em}\raisebox{-0.98ex}{\scalebox{2}{$\color{white}.$}}\hspace{-0.416em}\raisebox{+0.88ex}{$\color{white}.$}\hspace{0.46em}$}{}}

\nc{\bnu}{{\bar{ \nu}}}

\nc{\olO}{\bar{\calO}}

\nc{\al}{{\alpha}} 
\nc{\be}{{\beta}}
\nc{\ga}{{\gamma}} \nc{\Ga}{{\Gamma}}
 \nc{\hGa}{\hat{\Gamma}}
\nc{\ve}{{\varepsilon}} 
\nc{\la}{{\lambda}} \nc{\La}{{\Lambda}}
\nc{\om}{\omega} \nc{\Om}{\Omega} 
\nc{\sig}{{\sigma}} \nc{\Sig}{{\Sigma}}

\nc{\tnb}{\psi_{\rm tame}}
\nc{\oM}{\overline{{M}}}
\nc{\op}{{\on{op}}}
\nc{\ad}{{\on{ad}}}
\nc{\alg}{{\on{alg}}}
\nc{\Ad}{{\on{Ad}}}
\nc{\Adm}{{\on{Adm}}} \nc{\aff}{{\on{aff}}}
\nc{\Aut}{{\on{Aut}}}
\nc{\Bun}{{\on{Bun}}}
\nc{\cha}{{\on{char}}}
\nc{\der}{{\on{der}}}
\nc{\Der}{{\on{Der}}}
\nc{\diag}{{\on{diag}}}
\nc{\End}{{\on{End}}}
\nc{\Fl}{{\calF\!\ell}}
\nc{\Tr}{{\on{Transp}}}
\nc{\TR}{{\calT\!\calR}}
\nc{\Gal}{{\on{Gal}}}
\nc{\Gr}{{\on{Gr}}}
\nc{\rH}{{\on{H}}}
\nc{\Hom}{{\on{Hom}}}
\nc{\IC}{{\on{IC}}}
\nc{\id}{{\on{id}}}
\nc{\Id}{{\on{Id}}}
\nc{\ind}{{\on{ind}}}
\nc{\Ind}{{\on{Ind}}}
\nc{\Lie}{{\on{Lie}}}
\nc{\Pic}{{\on{Pic}}}
\nc{\pr}{{\on{pr}}}
\nc{\Res}{{\on{Res}}}
\nc{\res}{{\on{res}}} \nc{\Sat}{{\on{Sat}}}
\nc{\s}{{\on{sc}}}
\nc{\drv}{{\on{der}}}
\nc{\sgn}{{\on{sgn}}}
\nc{\Spec}{{\on{Spec}}}\nc{\Spf}{\on{Spf}} 
\nc{\Sph}{\on{Sph}}
\nc{\St}{{\on{St}}}
\nc{\tr}{{\on{tr}}}
\nc{\Mod}{{\mathrm{-Mod}}}
\nc{\Hilb}{{\on{Hilb}}} 
\nc{\Ext}{{\on{Ext}}} 
\nc{\vs}{{\on{Vec}}}
\nc{\ev}{{\on{ev}}}
\nc{\nO}{{\breve{\calO}}}
\nc{\tS}{{\tilde{S}}}
\nc{\spe}{{\on{sp}}}
\nc{\loc}{{\on{loc}}}
\nc{\Sym}{{\on{Sym}}}
\nc{\Cone}{{\on{C}}}
\nc{\syn}{{\on{syn}}}
\nc{\reg}{{\on{reg}}}
\nc{\colim}{{\on{colim}}}
\nc{\Norm}{{\on{N}}}

\nc{\nscrR}{{\mathscr{R}^{\on{nr}}}}

\nc{\GL}{{\on{GL}}}
\nc{\U}{{\on{U}}}
\nc{\Gl}{\on{Gl}} 
\nc{\GSp}{{\on{GSp}}}
\nc{\gl}{{\frakg\frakl}}
\nc{\SL}{{\on{SL}}} 
\nc{\SU}{{\on{SU}}} 
\nc{\SO}{{\on{SO}}}
\nc{\PGL}{{\on{PGL}}}

\nc{\Conv}{{\on{Conv}}}
\nc{\Rep}{{\on{Rep}}}
\nc{\Dom}{{\on{Dom}}}
\nc{\red}{{\on{red}}}
\nc{\act}{{\on{act}}}
\nc{\nr}{{\on{nr}}}
\nc{\ctf}{{\on{ctf}}}

\nc{\str}{{\on{-}}} 
\nc{\os}{{\bar{s}}}
\nc{\oeta}{{\bar{\eta}}}

\nc{\hookto}{\hookrightarrow}
\nc{\longto}{\longrightarrow}
\nc{\leftto}{\leftarrow}
\nc{\onto}{\twoheadrightarrow}
\nc{\lonto}{\twoheadleftarrow}
\newcommand*\isomto{%
  \renewcommand{\arraystretch}{0.1}
  \begin{array}[b]{c} {}_{\sim} \\ \longrightarrow \end{array}%
}

\nc{\uG}{{\underline{G}}}
\nc{\uA}{{\underline{A}}}
\nc{\uS}{{\underline{S}}}
\nc{\uT}{{\underline{T}}}
\nc{\uM}{{\underline{M}}}
\nc{\uP}{{\underline{P}}}
\nc{\uB}{{\underline{B}}}
\nc{\uN}{{\underline{N}}}

\nc{\ucG}{{\underline{\calG}}}
\nc{\ucA}{{\underline{\calA}}}
\nc{\ucS}{{\underline{\calS}}}
\nc{\ucT}{{\underline{\calT}}}
\nc{\ucalM}{{\underline{\calM}}}
\nc{\ucP}{{\underline{\calP}}}
\nc{\ucalN}{{\underline{\calN}}}

\nc{\bF}{{\breve{F}}}

\nc{\oFl}{{\overline{\Fl}}} 
\nc{\bU}{{\overline{U}}}
\nc{\tGr}{{\tilde{\Gr}}}
\nc{\cGr}{\calG\! r}
\nc{\oGr}{\overline{\on{Gr}}} 
\nc{\ocGr}{\overline{\calG\! r}}
\nc{\co}{{\colon}}
\nc{\sch}[1]{(Sch/{#1})}
\nc{\HypLoc}[1]{HypLoc({#1})}

\nc{\ohtimes}{\stackrel{!}{\otimes}}
\nc{\boxtilde}{\widetilde{\boxtimes}}
\nc{\vstar}{{\varhexstar}}

\nc{\Div}{\on{Div}}
\nc{\Sht}{\on{Sht}}
\nc{\Frob}{\on{Frob}}

\nc{\x}{\times}
\nc{\bsl}{\backslash}
\nc{\algQl}{{\bar{\bbQ}_\ell}}
\nc{\sF}{{\bar{F}}}
\nc{\nF}{{\breve{F}}}
\nc{\nW}{{W^{\on{nr}}}}
\nc{\sk}{{\bar{k}}}
\nc{\cont}{\on{c}}
\nc{\Supp}{\on{Supp}}
\nc{\blt}{\bullet}  
\nc{\dom}{\on{dom}}
\nc{\scon}{{\on{sc}}} 
\nc{\Affine}{\on{Aff}} 
\nc{\nscrA}{\mathscr{A}^{\on{nr}}} 
\nc{\nfraka}{{\bbf^{\on{nr}}}}
\nc{\ran}{{\rangle}}
\nc{\lan}{{\langle}}
\nc{\bk}{{\bar{k}}}
\nc{\tF}{{\tilde{F}}}
\nc{\sS}{{\bar{S}}}
\nc{\LG}{{^\text{L}\hspace{-0.04cm}G}}
\nc{\LL}{{^\text{L}\hspace{-0.07cm}L}}
\nc{\et}{{\text{\rm \'et}}}
\nc{\inv}{{\on{inv}}}
\nc{\Hecke}{{\on{Hecke}}}
\nc{\Isom}{{\on{Isom}}}
\nc{\oSht}{{\overline{\on{Sht}}}}
\nc{\umu}{{\underline \mu}}
\nc{\AIJ}{{\calO_X[{\scriptstyle{\calI\over \calJ}}]}}
\nc{\Proj}{{\on{Proj}}}
\nc{\Bl}{{\on{Bl}}}

\nc{\Pos}{{\on{Pos}}}
\nc{\Sets}{{\on{Sets}}}
\nc{\AffSch}{{\on{AffSch}}}
\nc{\Groups}{{\on{Groups}}}
\nc{\Gpds}{{\on{Groupoids}}}
\nc{\Sch}{{\on{Sch}}}
\nc{\fl}{{\on{flat}}}

\nc{\pot}[1]{ [\hspace{-0,5mm}[ {#1} ]\hspace{-0,5mm}] }
\nc{\rpot}[1]{ (\hspace{-0,7mm}( {#1} )\hspace{-0,7mm}) }

\nc{\defined}{\hspace{0.1cm}\stackrel{\text{\tiny \rm def}}{=}\hspace{0.1cm}}

\begin{document}

\title{N\'eron blowups and low-degree cohomological applications}
\author{by Arnaud Mayeux, Timo Richarz and Matthieu Romagny}

\address{Beijing International Center for Mathematical Research, Peking University, China}
\email{arnaud.mayeux@imj-prg.fr}
\address{Technical University of Darmstadt, Department of Mathematics, 64289 Darmstadt, Germany}
\email{richarz@mathematik.tu-darmstadt.de}
\address{Univ Rennes, CNRS, IRMAR - UMR 6625, F-35000 Rennes, France}
\email{matthieu.romagny@univ-rennes1.fr}


\maketitle

\begin{abstract} 
We define dilatations of general schemes and study their basic properties. 
Dilatations of group schemes are --in favorable cases-- again group schemes, called N\'eron blowups. 
We give two applications to their cohomology in degree zero (integral points) and degree one (torsors):
we prove a canonical Moy-Prasad isomorphism that identifies the graded pieces in the congruent filtration of $G$ with the graded pieces in its Lie algebra $\frakg$, and we show that many level structures on moduli stacks of $G$-bundles are encoded in torsors under N\'eron blowups of $G$.
\end{abstract}

\tableofcontents


\pagestyle{plain}

\section{Introduction}

\subsection{Motivation and goals} 
N\'eron blowups (or dilatations) provide a tool to modify group schemes over the fibers of a given Cartier divisor on the base. 
Classically, their integral points over discrete valuation rings appear as congruence subgroups of reductive groups over fields, see \cite[\S2.1.2]{Ana73}, \cite[p.~551]{WW80}, \cite[\S3.2]{BLR90}, \cite[\S\S7.2--7.4]{PY06} and \cite[\S 2.8]{Yu15}. 
Over two-dimensional base schemes N\'eron blowups also appear in \cite[p.~175]{PZ13}.

This note extends the theory of dilatations to general schemes which
--in the case of group schemes-- are also called N\'eron blowups. 
We give two applications involving their cohomology in degree zero (integral points) and in degree one (torsors). 
Our results concerning their integral points lead to a general form of an
isomorphism of Moy and Prasad, frequently used in representation theory.
Our results concerning their torsors show that these naturally encode many level structures 
on moduli stacks of bundles. 
This is used to construct integral models of moduli
stacks of shtukas with level structures as in \cite{Dri87} and
\cite{Var04} which for parahoric level structures
might be seen as function field analogues
of the integral models of Shimura varieties in \cite{KP18}.

\subsection{Results}
Let $S$ be a scheme. Let $S_0$ be an effective Cartier divisor on $S$,
i.e., a closed subscheme which is locally defined by a single non-zero
divisor. We denote by $\Sch_S^{S_0\text{-}\reg}$ the full subcategory
of schemes $T\to S$ such that $T|_{S_0}:=T\x_SS_0$ is an effective Cartier
divisor on $T$. This category contains all flat schemes over $S$.
For a group scheme $G\to S$ together with a closed subgroup
$H\subset G|_{S_0}$ over $S_0$, we define the contravariant functor
$\calG\co \Sch_S^{S_0\text{-}\reg}\to \Groups$ given by
all morphisms of $S$-schemes $T\to G$ such that the restriction
$T|_{S_0}\to G|_{S_0}$ factors through $H$.\medskip

\noindent{\bf Theorem.}
{\it \textup{(1)} The functor $\calG$ is representable by an open
subscheme of the full scheme-theoretic blowup of $G$ in $H$. 
The structure morphism $\calG\to S$ is an object in $\Sch_S^{S_0\text{-}\reg}$.
\textup{(\ref{blow.up.open.lemm}, \ref{blow.up.Cartier.lemm},
\ref{blow.up.rep.prop})} \smallskip\\
\textup{(2)} The canonical map $\calG\to G$ is affine. 
Its restriction over $S\bsl S_0$ induces an isomorphism
$\calG|_{S\bsl S_0}\cong G|_{S\bsl S_0}$. 
Its restriction over $S_0$ factors as $\calG|_{S_0}\to H\subset G|_{S_0}$.
\textup{(\ref{blow.up.Cartier.lemm}, \ref{Neron.blow.lemm})}\smallskip\\
\textup{(3)} If $H\to S_0$ has connected fibres and $H\subset G|_{S_0}$ is
regularly immersed, then $\calG|_{S_0}\to S_0$ has connected fibres.
\textup{(\ref{blow.up.smoothness.prop}, \ref{neron.blow.theo})}\smallskip\\
\textup{(4)} If $G\to S$, $H\to S_0$ are flat and \textup{(}locally\textup{)}
of finite presentation and $H\subset G|_{S_0}$ is
regularly immersed, then $\calG\to S$ is flat and \textup{(}locally\textup{)}
of finite presentation. If both $G\to S$,
$H\to S_0$ are smooth, then $\calG\to S$ is smooth.
\textup{(\ref{blow.up.smoothness.prop}, \ref{neron.blow.theo})}\smallskip\\
\textup{(5)} Assume that $\calG\to S$ is flat. Then its formation commutes
with base change $S'\to S$ in $\Sch_S^{S_0\text{-}\reg}$, and it carries
the structure of a group scheme such that the canonical map
$\calG\to G$ is a morphism of $S$-group schemes.
\textup{(\ref{blow.up.base.change.lemm}, \ref{neron.blow.theo})}\smallskip\\
\textup{(6)} Assume that $G\to S$ is flat, finitely presented and
$H\to S_0$ is flat, regularly immersed in $G|_{S_0}$. Locally over $S_0$,
there is an exact sequence of $S_0$-group schemes
$1 \to V \to \calG|_{S_0} \to H\to 1$ where $V$ is the vector bundle
given by restriction to the unit section of an explicit twist
of the normal bundle of $H$ in $G|_{S_0}$. If $H$ lifts to a flat
$S$-subgroup scheme of~$G$, this sequence is canonical; moreover
it exists globally and is canonically split.
\textup{(\ref{theo:dilatations of subgroups})}
}
\medskip

We call $\calG\to S$ the {\em N\'eron blowup \textup{(}{\em or}
dilatation\textup{)} of $G$ in $H$ along $S_0$}.
Note that $\calG\to S$ is a group object in $\Sch_S^{S_0\text{-}\reg}$ by (1),
but that $\calG\to S$ is a group scheme only if the self products $\calG\x_S\calG$ and $\calG\x_S\calG\x_S\calG$
are objects in $\Sch_S^{S_0\text{-}\reg}$ which holds for example in (5), cf.~\S\ref{neron.blow.def.sec} for details.
If $S$ is the spectrum of a discrete valuation ring and if $S_0$ is defined
by the vanishing of a uniformizer, then $\calG\to S$ is the group scheme
constructed in \cite[\S2.1.2]{Ana73}, \cite[p.~551]{WW80}, \cite[\S3.2]{BLR90}, \cite[\S\S7.2--7.4]{PY06} and \cite[\S 2.8]{Yu15}.
For an example of N\'eron blowups over two-dimensional
base schemes we refer to \cite[p.~175]{PZ13}, cf.~also Example \ref{neron.blow.exam}.

We point out that most of the foundations of the study of dilatations
can be settled in an absolute setting for schemes. That is, we initially
develop the theory of affine blowups (or dilatations) for closed
subschemes $Z$ in a scheme $X$ along a divisor $D$. It is only later
that we specialize to relative schemes (over some base $S$, with the
divisor coming from the base) and then further to group schemes.

The applications we give originate from a sheaf-theoretic viewpoint
on N\'eron blowups. Write $j\co S_0\hookto S$ the closed immersion of
the Cartier divisor, and assume that $G\to S$ and $H\to S_0$ are
flat, locally finitely presented groups. In this context, the dilatation
$\calG\to G$ sits in an exact sequence of sheaves of
pointed sets on the small syntomic site of $S$ (Lemma
\ref{lemma:short exact sequence}):
\[
\begin{tikzcd}[column sep=15]
1 \ar[r] & \calG \ar[r] & G \ar[r] & j_*(G_0/H) \ar[r] & 1,
\end{tikzcd}
\]
where $G_0:=G|_{S_0}$.
If $G\to S$ and $H\to S_0$ are smooth, then the sequence is
exact as a sequence of sheaves on the small \'etale site of $S$.
Considering the associated sequence on global sections,
we obtain the following theorem which generalizes and unifies several
results found in the literature (Remark \ref{literature.on.MP.remark}) under the name of {\em Moy-Prasad
isomorphisms}.\medskip

\noindent{\bf Corollary 1.}
{\it Let $r,s$ be integers such that $0\le r/2\le s\le r$.
Let $(\calO,\pi)$ be a henselian pair where $\pi\subset \calO$ is
an invertible ideal. Let $G$ be a smooth, separated
$\calO$-group scheme. Let $G_r$ be the $r$-th iterated dilatation of
the unit section and $\frakg_r$ its Lie algebra. If $\calO$ is local
or $G$ is affine, there is a canonical  isomorphism
$G_s(\calO)/G_r(\calO) \isomto \frakg_s(\calO)/\frakg_r(\calO)$.
\textup{(\ref{theoisomp})}
}\medskip

As another application, we are interested in comparing $\calG$-torsors
with $G$-torsors. In light of the above short exact sequence of
sheaves, there is an equivalence between the category of $\calG$-torsors
and the category of $G$-torsors equipped with a section of their
pushforward along $G\to j_*(G_0/H)$, see
\cite[Chap.~III, \S 3.2, Prop.~3.2.1]{Gi71}.

This has consequences for moduli of torsors over curves. We thus
specialize to the following setting, cf.~\S\ref{bundles.curves.sec}.
Assume that $X$ is a smooth, projective, geometrically irreducible curve
over a field $k$ with a Cartier divisor $N\subset X$, that $G\to X$ is a
smooth, affine group scheme and that $H\to N$ is a smooth closed subgroup
scheme of $G|_N$. In this case, the N\'eron blowup $\calG\to X$ is a smooth,
affine group scheme. Let $\Bun_G$ (resp.~$\Bun_\calG$) denote the moduli
stack of $G$-torsors (resp.~$\calG$-torsors) on $X$.  This is a quasi-separated,
smooth algebraic stack locally of finite type over $k$
(cf.~e.g.~\cite[Prop.~1]{He10} or \cite[Thm.~2.5]{AH19}).
Pushforward of torsors along $\calG\to G$ induces a morphism
$\Bun_\calG\to \Bun_G$, $\calE\mapsto \calE\x^\calG G$.
We also consider the stack $\Bun_{(G,H,N)}$ of $G$-torsors on $X$ with level-($H$,$N$)-structures, cf.~Definition \ref{neron.blow.groupoid.defi}.
Its $k$-points parametrize pairs $(\calE,\be)$ consisting of a $G$-torsor $\calE\to X$ and a section $\be$ of the fppf quotient $(\calE|_{N}/H)\to N$, i.e., $\be$ is a reduction of $\calE|_N$ to an $H$-torsor.\medskip

\noindent{\bf Corollary 2.} 
{\it There is an equivalence of $k$-stacks
\[
\Bun_\calG\;\overset{\cong}{\longto}\;\Bun_{(G,H,N)},\;\; \calE\mapsto (\calE\x^\calG G,\, \be_{\on{can}}),
\]
where $\be_{\on{can}}$ denotes the canonical reduction induced from the factorization $\calG|_N\to H\subset G|_N$ given by \textup{(2)} in the Theorem.}
(\ref{Neron.blow.equivalence.cor}, \ref{applications.bun.level.theo})\medskip

If $H=\{1\}$ is trivial, then $\Bun_{(G,H,N)}$ is the moduli stack of $G$-torsors equipped with level-$N$-structures. 
If $G\to X$ is reductive, if $N$ is reduced and if $H$ is a parabolic subgroup in $G|_N$, then $\Bun_{(G,H,N)}$ is the moduli stack of $G$-torsors with quasi-parabolic structures as in \cite{LS97}.
In this case the restriction of $\calG$ to the completed local rings of $X$ are parahoric group schemes in the sense of \cite{BT84} and the previous corollary was pointed out in \cite[\S 2.a.]{PR10}.
Thus, many level structures are encoded in torsors under N\'eron blowups.
This construction is also compatible with the adelic viewpoint, cf.~Corollary \ref{applications.Weil.uniformization.coro}.

Now assume that $k$ is a finite field. As a consequence of the
corollary one naturally obtains integral models for moduli stacks
of $G$-shtukas on $X$ with level structures over $N$ as
in \cite{Dri87} for $G=\GL_n$ and in \cite{Var04} for general split reductive $G$.
We thank Alexis Bouthier for informing us that this was already
pointed out in \cite[\S2.4]{NN08} for congruence level.
General properties of moduli stacks of shtukas for smooth, affine group schemes are studied in \cite{AH19}, \cite{AHab19} and \cite{Br}.
In \S\ref{applications.integral.models.sec} below, we make the connection between $G$-shtukas with level structures as in \cite{Dri87, Var04, Laf18} and $\calG$-shtukas as in \cite{AH19, AHab19, Br}. 
We expect the point of view of $\calG$-shtukas, as opposed to $G$-shtukas with level structures, to be fruitful for investigations also outside the case of parahoric level structures.

\ackn
We thank Anne-Marie Aubert,
Patrick Bieker, Paul Breutmann, Michel Brion, Colin Bushnell,
K\c{e}stutis \v{C}esnavi\v{c}ius, Laurent Charles, Cyril Demarche,
Antoine Ducros, Philippe Gille, Thomas Haines, Urs Hartl, Jochen Heinloth,
Eugen Hellmann, Laurent Moret-Bailly, Gopal Prasad,
 Beno\^it Stroh, Torsten Wedhorn and Jun Yu
for useful discussions around the subject of this note.
Also we thank Alexis Bouthier for pointing us to the reference \cite{NN08}.
\xackn

\section{Dilatations} 

In this section, we define dilatations and give some properties.
Dilatations (or affine blowups) are spectra of affine blowup algebras. We first introduce affine blowup algebras.

\subsection{Definition}\label{blow.up.definition.sec}
Fix a scheme $X$. Let $Z\subset D$ be closed subschemes in $X$,
and assume that $D$ is locally principal.  Denote by
$\calJ\subset \calI$ the associated quasi-coherent sheaves
of ideals in $\calO_X$ so that $Z=V(\calI)\subset V(\calJ)=D$. Let $\Bl_\calI\calO_X=\calO_X\oplus\calI\oplus\calI^2\oplus\ldots$ denotes the Rees algebra, it is a quasi-coherent $\bbZ_{\geq 0}$-graded $\calO_X$-algebra. If $\calJ=(b)$ is principal with $b\in \Ga(X,\calJ)$, then $\big(\Bl_\calI\calO_X\big)[\calJ^{-1}]:=\big(\Bl_\calI\calO_X\big)[b^{-1}]$ is well-defined independently of the choice of generators of $\calJ$.
If $\calJ$ is only locally principal, then we define the localization $\big(\Bl_\calI\calO_X\big)[\calJ^{-1}]$ by glueing. This $\calO_X$-algebra inherits a grading by giving local generators of $\calJ$ degree $1$. In other words, the graduation of $\big(\Bl_\calI\calO_X\big)[\calJ^{-1}]$ is given locally by $\on{deg}(i/b^k)= n-k$ for $i \in \calI ^n$ and $b$ a local generator of $\calJ$.

 \defi \label{blow.up.algebra.defi}
We use the following terminology:
\begin{enumerate} 
\item The {\it affine blowup algebra} of $\calO_X$ in $\calI$ along $\calJ$ is the quasi-coherent sheaf of $\calO_X$-algebras
\[
\AIJ\defined \Big[ \big(\Bl_\calI\calO_X\big)[\calJ^{-1}]\Big]_{\on{deg}=0},
\]
obtained as the subsheaf of degree $0$ elements in $\big(\Bl_\calI\calO_X\big)[\calJ^{-1}]$.
\item The {\it dilatation} (or {\it affine blowup}) of $X$ in $Z$ along $D$ is the $X$-affine scheme
\[
\Bl_Z^DX\defined \Spec\big(\AIJ\big).
\]
The subscheme $Z$, or the pair $(Z,D)$, is called the {\em center} of the dilatation.
\end{enumerate}
\xdefi

\rema
If $X$ is affine, affine blowup algebras are defined in \cite[\href{https://stacks.math.columbia.edu/tag/052P}{052P}]{StaPro}. In this case we denote $B:=\Ga(X,\calO_X)$, $I:=\Ga(X,\calI)$, $J:=\Ga(X,\calJ)$ and $\Bl_IB:=\Ga(X,\Bl_\calI\calO_X)=\oplus_{n\geq 0} I^n$. 
Moreover, if $J=(b)$ is principal, then $B[{I\over b}]:=\Ga(X,\AIJ)$ is the
algebra whose elements are equivalence classes of fractions $x/b^n$ with
$x\in I^n$, where two representatives $x/b^n$, $y/b^m$ with $x\in I^n$, $y\in I^m$
define the same element in $B[{I\over b}]$ if and only if there exists an
integer $l\geq 0$ such that
\begin{equation}\label{blow.up.equi.eq}
b^l(b^mx-b^ny)\;=\;0 \;\;\; \text{inside\; $B$.}
\end{equation}
 By \cite[\href{https://stacks.math.columbia.edu/tag/07Z3}{07Z3}]{StaPro},
\begin{align}
& \text{\rm the image of $b$ in $B[{\textstyle {I\over b}}]$ is a non-zero divisor,} \label{blow.up.nonzero.eq}\\
& bB[{\textstyle{I\over b}}]\;=\;IB[{\textstyle{I\over b}}], \;\;\text{and} \label{blow.up.divisor.eq}\\
& B[{\textstyle{I\over b}}][b^{-1}]\;=\; B[b^{-1}]. \label{blow.up.localization.eq}
\end{align}
In particular, the ring $B[{I\over b}]$ is the $B$-subalgebra of $B[b^{-1}]$ generated by fractions $x/b$ with $x\in I$.
\xrema

\subsection{Basic properties} 
We proceed with the notation from \S\ref{blow.up.definition.sec}. 
The following results generalize \cite[\S3.2, Prop.~1]{BLR90}.

\lemm \label{blow.up.open.lemm}
The affine blowup $\Bl_Z^DX$ is the open subscheme of the blowup $\Bl_ZX=\Proj(\Bl_\calI\calO_X)$ defined by the complement of $V_+(\calJ)$.
\xlemm
\pf
Our claim is Zariski local on $X$. 
We reduce to the case where $X=\Spec(B)$ is affine and $J=(b)$ is principal.
Then $B[{I\over b}]$ is the homogenous localization of $B\oplus I\oplus I^2\oplus\ldots$ at $b\in I$ viewed as an element in degree $1$, cf.~\cite[\href{https://stacks.math.columbia.edu/tag/052Q}{052Q}]{StaPro}.
This shows that $\Spec(B[{I\over b}])$ is the complement of $V_+(b)$ in $\Proj(\Bl_IB)$.
\xpf

\lemm\label{blow.up.Cartier.lemm}
As closed subschemes of $\Bl_Z^DX$, one has
\[
\Bl_Z^DX\x_XZ=\Bl_Z^DX\x_XD,
\]
which is an effective Cartier divisor on $\Bl_Z^DX$.
\xlemm
\pf
Our claim is Zariski local on $X$. 
We reduce to the case where $X=\Spec(B)$ is affine and $J=(b)$ is principal.
We have to show that $bB[{I\over b}]=IB[{I\over b}]$, and that $b$ is a non-zero divisor in $B[{I\over b}]$.
This is \eqref{blow.up.nonzero.eq} and \eqref{blow.up.divisor.eq} above.
\xpf

When $D$ is a Cartier divisor, we can also realize $\Bl_Z^DX$ as a
closed suscheme of the affine projecting cone. Recall that classically,
this cone is defined as the relative spectrum of the Rees algebra
$\Bl_{\calI}\calO_X$, so the blowup and its affine cone are
complementary to each other in the completed projective cone; see \cite[\S8.3]{EGA2}.
However, twisting the blowup algebras with the invertible ideal sheaf
$\calJ$ gives rise to different embeddings. Indeed, according to
\cite[\S8.1.1]{EGA2} we have a canonical isomorphic presentation
of the usual blowup as
$\Bl_ZX=\Proj(\oplus_{n\geq 0} \calI^n\otimes\calJ^{-n})$.
Here we define the {\em affine projecting cone} of $\Bl_ZX$
(with respect to the chosen presentation) as
\[
\Cone_ZX\defined\Spec\big(\oplus_{n\geq 0} \calI^n\otimes\calJ^{-n}\big).
\]

\lemm \label{blow.up.closed.in.cone.lemm}
If $D$ is a Cartier divisor, the affine blowup $\Bl_Z^DX$ is the closed
subscheme of the affine projecting cone $\Cone_ZX$ defined by the equation $\varrho-1$,
where $\varrho\in\calI\otimes\calJ^{-1}$ is the image of $1$ under the inclusion
$\calO_X=\calJ\otimes\calJ^{-1}\subset \calI\otimes\calJ^{-1}$.
\xlemm
\pf
Let $\scrA=\oplus_{n\geq 0} \calI^n\otimes\calJ^{-n}$.
There is a surjective morphism of sheaves of algebras $\scrA\to \AIJ$
defined by mapping a local section $i\otimes j^{-1}$ in degree 1 to $i/j$.
To check that $\varrho-1$ goes to zero and generates the kernel, we may
work locally on some affine open subscheme $U\subset X$ where the sheaf
$\calJ$ is generated by a section $b$. 
Let $t=b^\vee$ be the generator for $\calJ^{-1}$, dual to $b$.
Let $B=\Ga(X,\calO_X)$ and $I=\Ga(X,\calI)$.
Then the map $\scrA(U)\to \AIJ(U)$ is given by
\[
(\oplus_{n\geq 0} I^nt^n)\longto B[{\textstyle{I\over b}}],
\;\; \textstyle{\sum_{n\geq 0}} i_nt^n\mapsto \textstyle{\sum_{n\geq 0}} i_n/b^n.
\]
This induces an isomorphism
$(\oplus_{n\geq 0} I^nt^n)/(bt-1)\isomto B[{I\over b}]$.
\xpf

\subsection{Universal property}\label{blow.up.univ.ppty.section}
In this text we will use regular immersions in a possibly
non-noetherian setting where the reference \cite[\S16.9 and
\S19]{EGA4.4} is inadequate. In this case we refer to the Stacks
Project \cite{StaPro}. There, four notions of regularity are studied:
by decreasing order of strength, {\em regular},
{\em Koszul-regular}, {\em $H_1$-regular}, {\em quasi-regular} (see
Sections~\href{https://stacks.math.columbia.edu/tag/067M}{067M}
and \href{https://stacks.math.columbia.edu/tag/0638}{0638} in
{\em loc. cit.}). The useful ones for us are the first (which is
regularity in its classical meaning) and the third: an $H_1$-regular
sequence is a sequence whose Koszul complex has no homology in
degree~1. In \cite{StaPro} several results are stated under the
weakest $H_1$-regular assumption. For simplicity we will state our
results for regular immersions, although all of them hold also for
$H_1$-regular immersions.

Note that regularity and $H_1$-regularity coincide for sequences
composed of one element $x$, because for them the Koszul complex has
length one and the homology group in degree~1 is just the $x$-torsion. In particular, for locally principal subschemes the
three notions regular, Koszul-regular and $H_1$-regular are equivalent.

Let us denote by $\Sch_X^{D\text{-}\reg}$ the full subcategory of schemes $T\to X$ such that $T\x_XD \subset T$ is regularly immersed, or equivalently
is an effective Cartier divisor (possibly the empty set) on $T$.
If $T'\to T$ is flat and $T\to X$ is an object in this category, so is the composition $T'\to T\to X$. 
In particular, the category $\Sch_X^{D\text{-}\reg}$ can be equipped with the fpqc/fppf/\'etale/Zariski Grothendieck topology so that the notion of sheaves is well-defined.

As $\Bl_Z^DX\to X$ defines an object in $\Sch_X^{D\text{-}\reg}$ by Lemma \ref{blow.up.Cartier.lemm}, the contravariant functor
\begin{equation}\label{blow.up.represent.eq}
\Sch_X^{D\text{-}\reg}\to \Sets, \;\;\;(T\to X)\mapsto \Hom_{X\text{-Sch}}\big(T,\Bl_Z^DX\big)
\end{equation}
together with $\id_{\Bl_Z^DX}$ determines $\Bl_Z^DX\to X$ uniquely up to unique isomorphism.
The next proposition gives the universal property of dilatations.

\prop\label{blow.up.rep.prop}
The affine blowup $\Bl_Z^DX\to X$ represents the contravariant functor $\Sch_X^{D\text{-}\reg}\to \Sets$ given by
\begin{equation}\label{blow.up.iso.eq}
(f\co T\to X) \;\longmapsto\; \begin{cases}\{*\}, \; \text{if $f|_{T\x_XD}$ factors through $Z\subset X$;}\\ \varnothing,\;\text{else.}\end{cases}
\end{equation}
\xprop
\pf
Let $F$ be the functor defined by \eqref{blow.up.iso.eq}. 
If $T\to \Bl_Z^DX$ is a map of $X$-schemes, then the structure map $T\to X$ restricted to $T\x_XD$ factors through $Z\subset X$ by Lemma \ref{blow.up.Cartier.lemm}.
This defines a map
\begin{equation}\label{blow.up.map.sheaves.eq}
\Hom_{X\text{-Sch}}\big(\,\str\,,\Bl_Z^DX\big)\;\longto\; F
\end{equation}
of contravariant functors $\Sch_X^{D\text{-}\reg}\to \Sets$.
We want to show that \eqref{blow.up.map.sheaves.eq} is bijective when evaluated at an object $T\to X$ in $\Sch_X^{D\text{-}\reg}$.
As \eqref{blow.up.map.sheaves.eq} is a morphism of Zariski sheaves, we reduce to the case where both $X=\Spec(B)$, $T=\Spec(R)$ are affine and $J=(b)$ is principal.

For injectivity, let $g,g'\co B[{I\over b}]\to R$ be two $B$-algebra maps. 
We need to show $g=g'$. 
Indeed, since $B[b^{-1}]=B[{I\over b}][b^{-1}]$ by \eqref{blow.up.localization.eq}, we get $g[b^{-1}]=g'[b^{-1}]$.
As $b$ is a non-zero divisor in $R$ by assumption, this implies $g=g'$.
 
For surjectivity, consider an element in $F(\Spec(R))$ which corresponds to a ring morphism $g\co B\to R$ such that $I$ is contained in the kernel of $B\to R\to R/bR$.
We need to show that $g$ extends (necessarily unique) to an $B$-algebra morphism $\tilde g\co B[{I\over b}]\to R$.
Let $[x/b^n]$, $x\in I^n$ be a class in $B[{I\over b}]$.
Since $g(I^n)\subset (b^n)$ in $R$, the $b$-torsion freeness of $R$ implies that there is a unique element $r=r(x,n)\in R$ such that $g(x)=b^n\cdot r$.
We define $\tilde g([x/b^n]):=r(x,n)$.
This is well-defined: If $y/b^m$, $y\in I^m$ is another representative of $[x/b^n]$, then applying $g$ to equation \eqref{blow.up.equi.eq} yields $b^mg(x)=b^ng(y)$ in $R$.
It follows that $r(x,n)=r(y,m)$. 
Thus, $\tilde g$ is well-defined.
Similarly, one checks that $\tilde g$ defines a morphism of $B$-algebras.
\xpf

\subsection{Functoriality} \label{blow.up.base.functoriality.sec} 
Let $Z'\subset D' \subset X'$ be another triple as in \S\ref{blow.up.definition.sec}. 
A morphism $X'\to X$ such that its restriction to $D'$ (resp.~${Z'}$) factors through $D$ (resp.~$Z$) induces a unique morphism  $\Bl_{Z'}^{D'}X'\to \Bl_Z^DX$ such that the following diagram of schemes
\[
\begin{tikzpicture}[baseline=(current  bounding  box.center)]
\matrix(a)[matrix of math nodes, 
row sep=1.5em, column sep=2em, 
text height=1.5ex, text depth=0.45ex] 
{ 
\Bl_{Z'}^{D'}X'& \Bl_Z^DX \\ 
X'& X.\\}; 
\path[->](a-1-1) edge node[above] {} (a-1-2);
\path[->](a-2-1) edge node[above] {} (a-2-2);
\path[->](a-1-1) edge node[right] {} (a-2-1);
\path[->](a-1-2) edge node[right] {} (a-2-2);
\end{tikzpicture}
\] 
commutes. 
Indeed, the existence of $\Bl_{Z'}^{D'}X'\to \Bl_Z^DX$ follows directly from Definition \ref{blow.up.algebra.defi}.
The uniqueness can be tested Zariski locally on $X$ and $X'$ where it follows from \eqref{blow.up.nonzero.eq} and \eqref{blow.up.localization.eq}.

\subsection{Base change} \label{blow.up.base.change.sec}
Now let $X'\to X$ be a map of schemes, and denote by $Z'\subset D'\subset X'$ the preimage of $Z\subset D\subset X$.
Then $D'\subset X'$ is locally principal so that the affine blow $\Bl_{Z'}^{D'}X'\to X'$ is well-defined. 
By \S\ref{blow.up.base.functoriality.sec} there is a canonical morphism of $X'$-schemes
\begin{equation}\label{blow.up.base.change.eq}
\Bl_{Z'}^{D'}X'\;\longto\; \Bl_Z^DX\x_{X}X'.
\end{equation}

\lemm\label{blow.up.base.change.lemm}
If $\Bl_Z^DX\x_{X}X'\to X'$ is an object of $\Sch_{X'}^{D'\text{-}\reg}$, then \eqref{blow.up.base.change.eq} is an isomorphism.
\xlemm
\pf
Our claim is Zariski local on $X$ and $X'$.
We reduce to the case where both $X=\Spec(B)$, $X'=\Spec(B')$ are affine, and $J=(b)$ is principal.
We denote $Z'=\Spec(B'/I')$ and $D'=\Spec(B'/J')$. 
Then $J'=(b')$ is principal as well where $b'$ is the image of $b$ under $B\to B'$.
We need to show that the map of $B'$-algebras
\[
B'\otimes_BB\big[{\textstyle{I\over b}}\big]\;\longto\; B'\big[{\textstyle{I'\over {b'}}}\big]
\] 
is an isomorphism.
However, this map is surjective with kernel the $b'$-torsion elements
\cite[\href{https://stacks.math.columbia.edu/tag/0BIP}{0BIP}]{StaPro}.
As $b'$ is a non-zero divisor in $B'\otimes_BB\big[{{I\over b}}\big]$ by assumption, the lemma follows.
\xpf


\coro \label{blow.up.base.change.cor}
If the morphism $X'\to X$ is flat and has some property $\calP$ which is
stable under base change, then $\Bl_{Z'}^{D'}X'\to\Bl_{Z}^{D}X$ is flat
and has $\calP$.
\xcoro
\pf
Since flatness is stable under base change the projection $p\co\Bl_Z^DX\x_{X}X'\to \Bl_Z^DX$ is flat and has property $\calP$.
By Lemma \ref{blow.up.base.change.lemm}, it is enough to check that the closed subscheme $\Bl_Z^DX\x_XD'$ defines an effective Cartier divisor on $\Bl_Z^DX\x_{X}X'$.
But this closed subscheme is the preimage of the effective Cartier divisor $\Bl_Z^DX\x_XD$ under the flat map $p$, and hence is an effective Cartier divisor as well.
\xpf

\subsection{Exceptional divisor}
For closed subschemes $Z\subset D$ in $X$ with $D$ locally principal,
we saw in Lemma~\ref{blow.up.Cartier.lemm} that the preimage of the center
$\Bl_Z^DX\x_XZ=\Bl_Z^DX\x_XD$ is an effective Cartier divisor in $\Bl_Z^DX$.
It is called the {\em exceptional divisor} of the affine blowup.

In order to describe it, as before we denote by $\calI$ and $\calJ$
the sheaves of ideals of $Z$ and $D$ in $\calO_X$. Also we let
$\calC_{Z/D}=\calI/(\calI^2+\calJ)$ and $\calN_{Z/D}=\calC_{Z/D}^\vee$
be the conormal and normal sheaves of $Z$ in~$D$.

\prop \label{exceptional divisor}
Assume that $D\subset X$ is an effective Cartier divisor, and $Z\subset D$ is a regular immersion. Write $\calJ_Z:=\calJ|_Z$.
\begin{enumerate}
\item The exceptional divisor $\Bl_Z^DX\x_XZ\to Z$ is an affine space
fibration, Zariski locally over $Z$ isomorphic to
$\bbV(\calC_{Z/D}\otimes \calJ_Z^{-1})\to Z$.
\item If $H^1(Z,\calN_{Z/D}\otimes \calJ_Z)=0$ (for example if $Z$ is affine),
then $\Bl_Z^DX\x_XZ\to Z$ is {\em globally} isomorphic to
$\bbV(\calC_{Z/D}\otimes \calJ_Z^{-1})\to Z$.
\item If $Z$ is a transversal intersection in the sense that there
is a cartesian square of closed subschemes whose vertical maps are
regular immersions
\[
\begin{tikzcd}[column sep=30]
W \ar[r,hook] \arrow[rd, "\square",phantom] & X \\
Z \ar[r,hook] \ar[u,hook] & D \ar[u,hook]
\end{tikzcd}
\]
then $\Bl_Z^DX\x_XZ\to Z$ is {\em globally and canonically} isomorphic
to $\bbV(\calC_{Z/D}\otimes \calJ_Z^{-1})\to Z$.
\end{enumerate}
\xprop
\pf
Using Lemma~\ref{blow.up.closed.in.cone.lemm} we can view the affine blowup
$\Bl_Z^DX$ as the closed subscheme with equation $\varrho-1=0$ of the
affine projecting cone $\Cone_ZX=\Spec(\scrA)$ with
$\scrA=\oplus_{n\geq 0} \calI^n\otimes\calJ^{-n}$.
First we compute the preimage $\Cone_ZX\x_XZ$. We have:
\[
\scrA\otimes\calO_X/\calI
=\oplus_{k\ge 0}(\calI^k/\calI^{k+1})\otimes \calJ^{-k}.
\]
Since the immersions $Z\subset D\subset X$ are regular, the conormal sheaf
$\calC_{Z/X}=\calI/\calI^2$ is finite locally free, and the
canonical surjective morphism of $\calO_X$-algebras
\[
\Sym^\bullet(\calC_{Z/X})\longrightarrow
\text{Gr}_{\calI}^\bullet(\calO_X)
\]
is an isomorphism, that is,
$\Sym^k(\calC_{Z/X})\to \calI^k/\calI^{k+1}$ is an isomorphism
for each $k\ge 0$. Taking into account that $\calJ$ is an
invertible sheaf, we obtain
\[
\Sym^k(\calC_{Z/X}\otimes \calJ_Z^{-1})
\isomto (\calI^k/\calI^{k+1})\otimes \calJ^{-k}.
\]
It follows that
$\Spec(\scrA\otimes\calO_X/\calI)=\bbV(\calC_{Z/X}\otimes \calJ_Z^{-1})$
and that $\Bl_Z^DX\x_XZ$ is the closed subscheme cut out
by $\varrho-1$ inside the latter. Now consider the sequence
of conormal sheaves:
\[
0 \longto {\calC_{D/X}}|_Z \longto \calC_{Z/X} \longto \calC_{Z/D} \longto 0.
\]
By \cite[Prop.~16.9.13]{EGA4.4} or
\cite[\href{https://stacks.math.columbia.edu/tag/063N}{063N}]{StaPro}
in the non-noetherian case, this sequence
is exact and locally split (beware that $\calC_{Z/X}$ is denoted
$\mathscr{N}_{Z/X}$ in \cite{EGA4.4}). Using the fact
that $\calC_{D/X}|_Z\otimes \calJ_Z^{-1}\simeq\calO_Z$ is freely generated
by $\varrho$ as a subsheaf of $\calC_{Z/X}\otimes \calJ_Z^{-1}$,
we obtain an extension
\begin{equation}\label{exact.seq.conormal.eq}
0 \longto \varrho\,\calO_Z \longto \calC_{Z/X}\otimes \calJ_Z^{-1} \longto
\calC_{Z/D}\otimes \calJ_Z^{-1} \longto 0.
\end{equation}
Now we consider the three cases listed in the proposition.

\smallskip

\noindent (1) Since $\calC_{Z/D}\otimes \calJ_Z^{-1}$ is locally
free, locally over $Z$ we can choose a splitting of
the exact sequence \eqref{exact.seq.conormal.eq} of conormal sheaves:
\[
\calC_{Z/X}\otimes \calJ_Z^{-1}=\varrho\,\calO_Z\oplus \calC_{Z/D}\otimes \calJ_Z^{-1}.
\]
Mapping $\varrho\mapsto 1$ yields a morphism of $\calO_Z$-modules
\[
\calC_{Z/X}\otimes \calJ_Z^{-1} \longrightarrow
\calO_Z\oplus (\calC_{Z/D}\otimes \calJ_Z^{-1})
\subset \Sym(\calC_{Z/D}\otimes \calJ_Z^{-1})
\]
which extends to a surjection of algebras with
kernel $(\varrho-1)$:
\[
\Sym(\calC_{Z/X}\otimes \calJ_Z^{-1}) \longrightarrow
\Sym(\calC_{Z/D}\otimes \calJ_Z^{-1}).
\]
This identifies $\Bl_Z^DX\x_XZ$ with the affine space bundle
$\bbV(\calC_{Z/D}\otimes \calJ_Z^{-1})$, locally over $Z$.

\smallskip

\noindent (2) The exact sequence defines a class in
$\Ext^1_{\calO_X}(\calC_{Z/D}\otimes \calJ_Z^{-1},\calO_Z)$.
Because the conormal sheaf is locally free, we have:
\[
\Ext^1_{\calO_Z}(\calC_{Z/D}\otimes \calJ_Z^{-1},\calO_Z)
\simeq \Ext^1_{\calO_Z}(\calO_Z,\calC_{Z/D}^{\vee}\otimes \calJ_Z)
\simeq H^1(Z,\calN_{Z/D}\otimes \calJ_Z).
\]
By assumption this vanishes and we obtain a global splitting.
From this one concludes as before.

\smallskip

\noindent (3) If $Z$ is the transversal intersection of $W$ and $D$,
then we have two exact, locally split sequences:
\[
\begin{tikzcd}[row sep=2,column sep=15]
0 \ar[rd] & & & & 0 \\
& \calC_{D/X}|_Z \ar[rd] & & \calC_{Z/W} \ar[ru] & \\
& & \calC_{Z/X} \ar[rd] \ar[ru] & & \\
& \calC_{W/X}|_Z \ar[ru] \ar[rr,dashed] & & \calC_{Z/D} \ar[rd] & \\
0 \ar[ru] & & & & 0.
\end{tikzcd}
\]
We claim that the dashed arrow is an isomorphism. To see this,
write $\calI,\calJ,\calK$ the defining ideals of $Z,D,W$. The composition
$\calC_{W/X}|_Z\to \calC_{Z/D}$ is the following map:
\[
\calK/(\calK^2+\calI\calK) \longrightarrow \calI/(\calI^2+\calJ).
\]
From the fact that $\calI=\calJ+\calK$ we deduce:
\begin{itemize}
\item $\calK^2+\calI\calK=\calK^2+\calJ\calK$, hence
$\calK/(\calK^2+\calI\calK)=\calK/(\calK^2+\calJ\calK)$.
\item $\calI^2=\calJ^2+\calJ\calK+\calK^2$, hence $\calI^2+\calJ=\calK^2+\calJ$
and $\calI/(\calI^2+\calJ)=(\calJ+\calK)/(\calK^2+\calJ)=\calK/(\calK^2+\calJ\cap \calK)$.
\end{itemize}
Hence, the map above is an isomorphism if and only if $\calJ\calK=\calJ\cap \calK$,
which holds because $W$ cuts $D$ transversally (this is another way of saying
that a local equation for $\calJ$ remains
a non-zero divisor in $\calO_W$).
This provides a canonical splitting
$\calC_{Z/X}=\calC_{D/X}|_Z\oplus \calC_{Z/D}$.
One concludes as before.
\xpf

\rema
In the course of the proof, we saw that the exceptional divisor
has the following explicit description: as an affine space fibration
over $Z$, its local sections over an open $U\subset Z$ are the
$\calO_U$-linear maps $\varphi:\calC_{U/X}\otimes \calJ_U^{-1}\to\calO_U$
such that $\varphi(\varrho)=1$.
\xrema

\subsection{Iterated dilatations}
Here we study the behaviour of dilatations under iteration.
Namely, we will
prove that when the center $Z$ of the affine blowup is a transversal
intersection $W\cap D$, it can be dilated any finite number of times and
the result of~$r$ dilatations can be seen as the dilatation
of the single ``thickened'' center $rZ$ (to be defined below)
inside the multiple Cartier divisor $rD$. To make this precise,
we first study the lifting of subschemes along an affine blowup.

\lemm \label{lemma:lifting immersions}
Let $Z\subset D\subset X$ be closed subschemes
with $D\subset X$ a Cartier divisor.
Let $\iota\co W\hookto X$ be an immersion such that $\iota^{-1}(Z)=\iota^{-1}(D)$
is a Cartier divisor in $W$. Set $X':=\Bl_Z^DX$ and let $\iota'\co W\to X'$
be the lift of $\iota$ given by the universal property of the dilatation.
\begin{enumerate}
\item If $\iota$ is an open immersion, then $\iota'$ is an
open immersion.
\item If $\iota$ is a closed immersion, then $\iota'$ is a
closed immersion.
\item Write $\iota$ as the composition
$\begin{tikzcd}[column sep=6mm]
W \ar[r,hook,closed] & U \ar[r,hook,open] & X
\end{tikzcd}$
where $U\subset X$ is the largest open subscheme such that $W$ is a
closed subscheme of $U$. Let $\calJ$ \textup{(}resp. $\calK$\textup{)}
be the ideal sheaf of $D$ \textup{(}resp. $W$\textup{)} in $U$.
Then $\iota'$ is the composition
$\begin{tikzcd}[column sep=6mm]
W \ar[r,hook,closed] & U' \ar[r,hook,open] & X'
\end{tikzcd}$
where $U'=X'\times_XU$ is the preimage of $U$ and $W\hookto U'$
is a closed immersion
with sheaf of ideals $\calK\calO_{U'}\otimes (\calJ\calO_{U'})^{-1}$.
\end{enumerate}
\xlemm

\pf
(1) In this case $\iota$ is flat, and the formation of the dilatation
commutes with base change. That is, the canonical morphism of
$W$-schemes
\[
\Bl_{Z\cap W}^{D\cap W}W \longto X'\times_XW
\]
is an isomorphism. But by the assumption
$\Bl_{Z\cap W}^{D\cap W}W\to W$ is an isomorphism, and this
identifies $W\to X'$ with the preimage of $W\hookto X$ in $X'$.

\smallskip

\noindent (2) Let $\calK\subset \calO_X$ be the ideal sheaf of $W$.
We will prove that $\iota'\co W\to X'$ is a closed immersion with ideal
sheaf $\calK\calO_{X'}\otimes (\calJ\calO_{X'})^{-1}$. First of all
$\iota'$ is automatically a monomorphism of schemes, and a proper map
because $\iota$ is proper
and $X'\to X$ is separated. Therefore $\iota'$ is a closed immersion by
\cite[Cor.~18.12.6]{EGA4.4}. The computation of the ideal sheaf is a
local matter so we can suppose that $X=\Spec(A)$ is affine and
the ideal sheaf $\calJ$ is generated by a section $b$. We write
$I,J,K\subset A$ the ideals defining $Z,D,W$ and $t:=b^\vee$ the
generator of $\calJ^{-1}$, dual to $b$. The assumptions of the lemma mean that
$I+K=J+K$ and $b$ is a non-zero divisor in $A/K$. From
Lemma~\ref{blow.up.closed.in.cone.lemm}, we know
that $X'$ is the spectrum of the ring
\[
A'=(\oplus_{e\ge 0} I^et^e)/(bt-1).
\]
In the present local situation, the map $\iota'\co W\to X'$ is
given by a lifting of
$\iota^\sharp\co A\to A/K$ to a map $(\iota')^\sharp\co A'\to A/K$.
Since $A'$ is generated by $It$ as an $A$-algebra, this map
is determined by the formula $(\iota')^\sharp(it)=j^\sharp(a)$,
for all $i\in I$ written $i=ab+k\in I\subset bA+K$.
In particular we see that $(\iota')^\sharp(Kt)=0$. Now working
modulo $(bt-1)+Kt$ in the ring $C=\oplus_{e\ge 0} I^et^e$, we have:
\[
It\subset bt A+Kt\equiv A+Kt \equiv A
\]
which sits in the degree~0 part of $C$, whence a surjection
$A\hookto C\to A'/KtA'$.
Moreover $bKt \equiv K$ implies that $K$ in degree~0
belongs to the ideal generated by $Kt$, hence finally
$A'/KtA'\isomto A/K$ as desired.

\smallskip

\noindent (3) This is the conjunction of (1) and (2).
\xpf

In view of the preceding lemma, if we fix a closed subscheme
$i\co W\hookto X$ such that $i^{-1}(Z)=i^{-1}(D)$ is a Cartier divisor
in $W$, we will be able to lift $W$ to the dilatation and hence
iterate the process. So we place ourselves in the following situation.

\begin{assum} \label{assumption2}
The schemes $Z\subset D\subset X$ sit in a cartesian diagram
of closed subschemes
\[
(\mathscr{D}): \qquad
\begin{tikzcd}[column sep=30]
W \ar[r,hook,"i"] \arrow[rd, "\square",phantom] & X \\
Z \ar[r,hook] \ar[u,hook] & D \ar[u,hook]
\end{tikzcd}
\]
such that the vertical maps are Cartier divisor inclusions.
\end{assum}

In this situation we can construct a sequence of dilatations
\[
\dots \longto X_r \longto X_{r-1} \longto \dots \dots
\longto X_1 \longto X_0=X
\]
and closed immersions $i_r\co W\hookto X_r$, as follows. We let
$\mathscr{D}_0 = \mathscr{D}$,
$X_0=X$, $D_0=D$, and $i_0=i\co W\hookto X_0$. Let $u_1\co X_1\to X_0$
be the dilatation of $Z$ in $(X_0,D_0)$, and $D_1$ the preimage
of $D_0$ in $X_1$.
Since $\mathscr{D}_0$ is a cartesian diagram of closed subschemes,
we have $i_0^{-1}(Z)=i_0^{-1}(D)=X$ which is a Cartier divisor in $W$.
So by the universal property
and Lemma~\ref{lemma:lifting immersions}, there is a closed
immersion $i_1\co W\to X_1$ lifting $i_0$. Moreover,
\[
(i_1)^{-1}(D_1)=(i_1)^{-1}(u_1^{-1}(D))=i^{-1}(D)=Z.
\]
That is, we again have a cartesian diagram
\[
(\mathscr{D}_1): \qquad
\begin{tikzcd}[column sep=30]
W \ar[r,hook,"i_1"] \arrow[rd, "\square",phantom] & X_1 \\
Z \ar[r,hook] \ar[u,hook] & D_1 \ar[u,hook]
\end{tikzcd}
\]
where the vertical maps are Cartier divisor inclusions.
Our sequence is obtained by iterating this construction.

\lemm \label{lemma:1}
Under Assumption~\ref{assumption2}, denote by $\calI,\calJ,\calK$
the ideal sheaves of $Z,D,W$ in $\calO_X$. Let $rD$ be the $r$-th
multiple of $D$ as a Cartier divisor, and $rZ:=W\cap rD$.
Then the dilatation $v_r\co X'_r\to X$ of $(rZ,rD)$ in $X$ is
characterized as being universal among all morphisms
$V\to X$ with the following two properties:
\begin{enumerate}
\item[\rm (i)] $\calJ\calO_V$ is an invertible sheaf,
\item[\rm (ii)] $\calK\calO_V$ is divisible by
$\calJ^r\calO_V$, that is we have $\calK\calO_V=\calJ^r\calO_V\cdot \calK_r$
for some sheaf of ideals $\calK_r\subset\calO_V$.
\end{enumerate}
\xlemm

\pf
The defining properties of the dilatation~$v_r$ say that it is
universal among morphisms $V\to X$ such that
$rZ\times_XV=rD\times_XV$ is a Cartier divisor. Since the ideal
sheaves of $rD$ and $rZ$ are $\calJ^r$ and $\calJ^r+\calK$ respectively,
these properties mean that the ideal $\calJ^r\calO_V$ is invertible
and $\calJ^r\calO_V=(\calJ^r+\calK)\calO_V$. But the properties
``$\calJ$ is invertible'' and ``$\calJ^r$ is invertible'' are equivalent,
as follows from the isomorphism between the blow-up of $\calJ$ and
the blow-up of $\calJ^r$, see \cite[Def.~8.1.3]{EGA2}. This takes
care of~(i). Besides, $\calJ^r\calO_V=(\calJ^r+\calK)\calO_V$ means that
$\calK\calO_V\subset \calJ^r\calO_V$ and in the situation where
$\calJ\calO_V$ is invertible, this is the same as saying that
\[
\calK\calO_V=\calJ^r\calO_V \cdot \calK_r
\]
with $\calK_r=(\calK\calO_V:\calJ^r\calO_V)$ as an ideal of $\calO_V$.
(Note that in this case
$(\calK\calO_V:\calJ^r\calO_V)\simeq (\calK\otimes \calJ^{-r})\calO_V$ as
an $\calO_V$-module.) This takes care of~(ii).
\xpf

\prop \label{prop:dilatation of thick subscheme}
In the situation of Assumption~\ref{assumption2}, let
\[
\dots \longto X_r \longto X_{r-1} \longto \dots \dots
\longto X_1 \longto X_0=X
\]
be the sequence of dilatations constructed above.
Let $rD$ be the $r$-th multiple of $D$ as a Cartier divisor,
and $rZ:=W\cap rD$.
Then the composition $X_r\to X$ is the dilatation of
$(rZ,rD)$ inside $X$.
\xprop

\pf
According to Lemma~\ref{lemma:1}, dilating $(rZ,rD)$ means
making $\calJ$ invertible and $\calK$ divisible by $\calJ^r$, all of this
in a universal way. This can be done by steps:
\begin{itemize}
\item make $\calJ$ invertible and make $\calK_0=\calK$ divisible by $\calJ$;
\item keep $\calJ$ invertible and make
$\calK_1:=(\calK_0:\calJ)\simeq \calK\otimes \calJ^{-1}$ divisible by $\calJ$;
\item keep $\calJ$ invertible and make
$\calK_2:=(\calK_1:\calJ)\simeq \calK_1\otimes \calJ^{-1}\simeq
\calK\otimes \calJ^{-2}$
divisible by $\calJ$; etc, and finally
\item keep $\calJ$ invertible and make
$\calK_{r-1}:=(\calK_{r-2}:\calJ)\simeq K\otimes \calJ^{-(r-1)}$
divisible by $\calJ$.
\end{itemize}
In view of Lemma~\ref{lemma:lifting immersions}, these steps amount to:
\begin{itemize}
\item dilate $Z$ in $(X,D)$,
\item dilate $Z$ in $(X_1,D_1)$,
\item dilate $Z$ in $(X_2,D_2)$, etc.~until
\item dilate $Z$ in $(X_{r-1},D_{r-1})$.
\end{itemize}
In this way we see the equivalence between the dilatation
of the thick pair $(rZ,rD)$ and the sequence of dilatations
of $Z$ constructed after~\ref{assumption2}.
\xpf

\subsection{Flatness and smoothness}
Flatness and smoothness properties of blowups are discussed in
\cite[\S 19.4]{EGA4.4}. 
Here we need slightly different versions.
We proceed with the notation from \S\ref{blow.up.definition.sec}.
We assume further that there exists a scheme $S$ under $X$ together with a locally principal closed subscheme $S_0\subset S$ fitting into a commutative diagram of schemes
\begin{equation}\label{blow.up.smoothness.diag}
\begin{tikzpicture}[baseline=(current  bounding  box.center)]
\matrix(a)[matrix of math nodes, 
row sep=1.5em, column sep=2em, 
text height=1.5ex, text depth=0.45ex] 
{ 
Z& D& X \\ 
& S_0 & S,\\}; 
\path[->](a-1-1) edge node[above] {} (a-1-2);
\path[->](a-1-2) edge node[above] {} (a-1-3);
\path[->](a-2-2) edge node[above] {} (a-2-3);
\path[->](a-1-1) edge node[right] {} (a-2-2);
\path[->](a-1-2) edge node[right] {} (a-2-2);
\path[->](a-1-3) edge node[right] {} (a-2-3);
\end{tikzpicture}
\end{equation}
where the square is cartesian, that is $D\to X_0:=X\times_SS_0$
is an isomorphism.

\lemm \label{flatness.by.local.criterion.lemma}
Assume that $S_0$ is an effective Cartier divisor in $S$.
Let $f:Y\to S$ be a morphism of schemes such that $Y_0:=Y\times_SS_0$
is a Cartier divisor in $Y$. Assume that both restrictions of $f$
above $S\bsl S_0$ and $S_0$ are flat. If one of the following holds:
\begin{itemize}
\item[(i)] $S,Y$ are locally noetherian,
\item[(ii)] $Y\to S$ is locally of finite presentation,
\end{itemize}
then $f$ is flat.
\xlemm

\pf
Since by assumption $u$ is flat at all points above the open
subscheme $S\bsl S_0$, it is enough to prove that $u$ is flat at
all points $y\in Y$ lying above a point $s\in S_0$.

In case (i), the local criterion for flatness
\cite[Chap.~O, 10.2.2]{EGA3.1} (cf.~also \cite[\href{https://stacks.math.columbia.edu/tag/00ML}{00ML}]{StaPro}) shows that $\calO_{S,s}\to\calO_{Y,y}$ is flat
and we are done.

In case (ii), we may localize around $y$ and $s$ and hence assume
that $Y$ and $S$ are affine and small enough so that the ideal
sheaf of $S_0$ in $S$ is generated by an element $f\in A=\Ga(S,\calO_S)$.
We write $A=\colim A_i$ as the union of its subrings of finite type
over $\bbZ$. In each $A_i$, the element $f$ is a non-zero divisor.
Write $S_i:=\Spec(A_i)$ and $S_{i,0}:=\Spec(A_i/f)$. Using the results
of \cite[Chap.~IV, \S8]{EGA4.3} on limits of schemes, we can find
an index $i$ and a morphism of finite presentation $Y_i\to S_i$ such that
$Y_{i,0}:=Y_i\times_{S_i} S_{i,0}$ is a Cartier divisor in $Y_i$
and $Y_{i,0}\to S_{i,0}$ is flat, and such that the situation
$(S,S_0,Y)$ is a pullback of $(S_i,S_{i,0},Y_i)$ by $S\to S_i$.
More in detail, using the following
results we find indices which we increase at each step in order to
have all the conditions simultaneously met: use \cite[Thm.~8.8.2]{EGA4.3}
to find morphisms $Y_i\to S_i$ and $Y_{i,0}\to S_{i,0}$, use
\cite[Cor.~8.8.2.5]{EGA4.3} to make $Y_i\times_{S_i} S_{i,0}$
and $Y_{i,0}$ isomorphic over $S_{i,0}$, use
\cite[Thm.~11.2.6]{EGA4.3} to ensure $Y_{i,0}\to S_{i,0}$ flat,
and use \cite[Prop.~8.5.6]{EGA4.3} to ensure that $f$ is a
non-zero divisor in $\calO_{Y_i}$, i.e. $Y_{i,0}\subset Y_i$ is a Cartier
divisor. Since $A_i$ is noetherian, for $(S_i,S_{i,0},Y_i)$ we can apply
case~(i) and the result follows by base change.
\xpf

For our conventions on regular immersions, the reader is referred
back to \S\ref{blow.up.univ.ppty.section}.

\prop \label{blow.up.smoothness.prop}
Assume that $S_0$ is an effective Cartier divisor on $S$.
\begin{enumerate}
\item If $Z\subset D$ is regular, then
$\Bl_Z^DX\to X$ is of finite presentation.
\item If $Z\subset D$ is regular, the fibers of
$\Bl_Z^DX\times_SS_0\to S_0$ are connected \textup{(}resp.
irreducible, geometrically connected, geometrically irreducible\textup{)}
if and only if the fibers of $Z\to S_0$ are.
\item If $X\to S$ is flat and if moreover one of the following holds:
\begin{itemize}
\item[(i)] $Z\subset D$ is regular, $Z\to S_0$ is flat
and $S,X$ are locally noetherian,
\item[(ii)] $Z\subset D$ is regular, $Z\to S_0$ is flat and
$X\to S$ is locally of finite presentation,
\item[(iii)] the local rings of $S$ are valuation rings,
\end{itemize}
then $\Bl_Z^DX\to S$ is flat.
\item If both $X\to S$, $Z\to S_0$ are smooth, then $\Bl_Z^DX\to S$ is smooth.
\end{enumerate}
\xprop

\pf
For (1) recall that the blowup of a regularly immersed subscheme has
an explicit structure, where generating relations between local generators
of the blown up ideal are the obvious ones, in finite number; see
\cite[\href{https://stacks.math.columbia.edu/tag/0BIQ}{0BIQ}]{StaPro}.
This shows that $\Bl_Z^DX\to X$ is locally of finite presentation.
Being also affine, it is of finite presentation.

For (2) about connectedness and irreducibility, recall from
Proposition~\ref{exceptional divisor}~(1) that the exceptional divisor
is an affine space fibration over $Z$. In particular it is a
submersion, so that the elementary topological lemma
\cite[Lem.~4.4.2]{EGA4.2} (cf.~also \cite[\href{https://stacks.math.columbia.edu/tag/0377}{0377}]{StaPro})
gives the assertion.

For (3)(i)-(ii), we apply Lemma~\ref{flatness.by.local.criterion.lemma}
to $Y:=\Bl_Z^DX$. The preimage of $S_0$ under the affine blowup
$f\co Y\to S$ is equal to $\Bl_Z^DX\x_XD=\Bl_Z^DX\x_XZ$ by
Lemma \ref{blow.up.Cartier.lemm}.
This implies that the restriction $f|_{f^{-1}(S\bsl S_0)}$ is equal
to $X\bsl D\to S\bsl S_0$ which is flat by assumption. It remains
to show flatness in points of $\Bl_Z^DX$ lying over $S_0$.
For this note that the restriction $f|_{f^{-1}(S_0)}$ factors as 
\begin{equation}\label{blow.up.smoothness.prop.eq1}
\Bl_Z^DX\x_XD=\Bl_Z^DX\x_XZ\longto Z\longto S_0,
\end{equation}
where the first map is smooth by Proposition~\ref{exceptional divisor}
and the second map is flat by assumption. Then
Lemma~\ref{flatness.by.local.criterion.lemma} applies and gives
flatness of $Y\to S$.

For (3)(iii) we can work locally at a point of $S$ and hence assume
that $S$ is the spectrum of a valuation ring $R$. We use the fact that
flat $R$-modules are the same as torsionfree $R$-modules. Locally over
an open subscheme $\Spec(B)\subset X$, the Rees algebra $\Bl_IB=B[It]$
is a subalgebra of the polynomial algebra $B[t]$ and the affine blowup
algebra is a localization of the latter. It follows that if $B$ is
$R$-torsionfree then the affine blowup algebra also, hence it is flat.

For (4) assume that $X\to S$ and $Z\to S_0$ are smooth.
Then (4) follows from
 \cite[Thm.~17.5.1]{EGA4.4} (cf.~also \cite[\href{https://stacks.math.columbia.edu/tag/01V8}{01V8}]{StaPro})
once we know that $\Bl_Z^DX\to S$ is locally of finite presentation,
flat and has smooth fibers.
Applying \cite[Prop.~19.2.4]{EGA4.4} to the commutative triangle in \eqref{blow.up.smoothness.diag} we see that $Z\subset D$ is regularly
immersed. Therefore, $\Bl_Z^DX\to S$ is flat and locally of
finite presentation by parts (1) and (3).
The smoothness of the fiber over points in $S\bsl S_0$ is clear, and follows from \eqref{blow.up.smoothness.prop.eq1} over points in $S_0$. 
This proves (4).
\xpf

\section{N\'eron blowups}
We extend the theory of N\'eron blowups of affine group schemes over discrete valuation rings as in \cite[2.1.2]{Ana73}, \cite[p.~551]{WW80}, \cite[\S3.2]{BLR90}, \cite[\S 2.8]{Yu15} and \cite[\S\S7.2--7.4]{PY06} to group schemes over arbitrary bases.

\subsection{Definition} \label{neron.blow.def.sec}
Let $S$ be a scheme, and let $G\to S$ be a group scheme.
Let $S_0\subset S$ be a locally principal closed subscheme,
and consider the base change $G_0:=G\x_SS_0$. 
Let $H\subset G_0$ be a closed subgroup scheme over $S_0$. 
Let $\calG:=\Bl_H^{G_0}G\to G$ be the dilatation of $G$ in $H$ along the locally principal, closed subscheme $G_0\subset G$ in the sense of Definition \ref{blow.up.algebra.defi}.
In this case, we also call $\calG\to S$ the {\it N\'eron blowup of $G$ in $H$ \textup{(}along $S_0$\textup{)}}.
We denote by $\calG_0:=\calG\x_SS_0\to S_0$ its exceptional divisor.

Let $\Sch_S^{S_0\text{-}\reg}$ be the full subcategory of schemes $T\to S$
such that $T_0:=T\x_SS_0$ defines an effective Cartier divisor on $T$.
By Lemma \ref{blow.up.Cartier.lemm} the structure morphism $\calG\to S$
defines an object in $\Sch_S^{S_0\text{-}\reg}$.

\lemm\label{Neron.blow.lemm} 
Let $\calG\to S$ be the N\'eron blowup of $G$ in $H$ along $S_0$.
\begin{enumerate}
\item The scheme $\calG\to S$ represents the contravariant functor
$\Sch_S^{S_0\text{-}\reg}\to \Sets$ given for $T\to S$ by the set of all $S$-morphisms $T\to G$ such that the induced morphism $T_0\to G_0$ factors through $H\subset G_0$.
\item The map $\calG\to G$ is affine. 
Its restriction over $S\bsl S_0$ induces an isomorphism $\calG|_{S\bsl S_0}\cong G|_{S\bsl S_0}$.
Its restriction over $S_0$ factors as $\calG_0\to H\subset  G_0$.
\end{enumerate}
\xlemm
\pf
Part (1) is a reformulation of Proposition \ref{blow.up.rep.prop}, and (2) is immediate from Lemmas \ref{blow.up.open.lemm} and \ref{blow.up.Cartier.lemm}.
\xpf

By virtue of Lemma \ref{Neron.blow.lemm} (1) the (forgetful) map
$\calG\to G$ defines a subgroup functor when restricted to the
category $\Sch_S^{S_0\text{-}\reg}$. As $\calG\to S$ is an object in
$\Sch_S^{S_0\text{-}\reg}$, it is a group object in this category. 
Here we note that products in the category $\Sch_S^{S_0\text{-}\reg}$
exist and are computed as $\Bl_{S_0}(X_1\x_SX_2)$ by the universal property of the blowup \cite[\href{https://stacks.math.columbia.edu/tag/085U}{085U}]{StaPro}. 
This is the closed subscheme of $X_1\x_SX_2$ which is locally defined by the ideal of $a$-torsion elements for a local equation $a$ of $S_0$ in $S$. 
In particular, if $\calG\to S$ is flat, then it is equipped with the structure of a group scheme such that $\calG\to G$ is a morphism of $S$-group schemes.

\subsection{Properties} \label{neron.blow.prop.sec}
We continue with the notation of \S\ref{neron.blow.def.sec}. 
Additionally assume that $S_0$ is an effective Cartier divisor in $S$.
Again recall our conventions on regular immersions from
\S\ref{blow.up.univ.ppty.section}. The following summarizes the
main properties of N\'eron blowups.

\theo\label{neron.blow.theo}
Let $\calG\to G$ be the N\'eron blowup of $G$ in $H$ along $S_0$.
\begin{enumerate}
\item If $G\to S$ is \textup{(}quasi-\textup{)}affine, then
$\calG\to S$ is \textup{(}quasi-\textup{)}affine.
\item If $G\to S$ is \textup{(}locally\textup{)} of finite
presentation and $H\subset G_0$ is regular, then
$\calG\to S$ is \textup{(}locally\textup{)} of finite presentation.
\item If $H\to S_0$ has connected fibres and $H\subset G_0$ is
regular, then $\calG\x_SS_0\to S_0$ has connected fibres.
\item Assume that $G\to S$ is flat and one of the following holds:
\begin{itemize}
\item[(i)] $H\subset G_0$ is regular, $H\to S_0$ is flat
and $S,G$ are locally noetherian,
\item[(ii)] $H\subset G_0$ is regular, $H\to S_0$ is flat and $G\to S$
is locally of finite presentation,
\item[(iii)] the local rings of $S$ are valuation rings,
\end{itemize}
then $\calG\to S$ is flat.
\item If both $G\to S$, $H\to S_0$ are smooth, then $\calG\to S$ is smooth.
\item Assume that $\calG\to S$ is flat. If $S'\to S$ is a scheme such
that $S_0':=S'\x_SS_0$ is an effective Cartier divisor on $S'$, then the base change $\calG\x_SS'\to S'$ is the N\'eron blowup of $G\x_SS'$ in $H\x_{S_0}S_0'$ along $S_0'$.
\end{enumerate}
In cases \textup{(4)} and \textup{(5)},
the map $\calG\to S$ is a group scheme.
\xtheo
\pf
The map $\calG\to G$ is affine by Lemma \ref{Neron.blow.lemm}~(2)
which implies (1). Items (2) to (5) are a direct transcription of Proposition~\ref{blow.up.smoothness.prop}, noting for (3) that
for schemes equipped with a section the properties ``with connected
fibers'' and ``with geometrically connected fibers'' are equivalent
\cite[Cor.~4.5.14]{EGA4.2} (cf.~also \cite[\href{https://stacks.math.columbia.edu/tag/04KV}{04KV}]{StaPro}).  
Part (6) follows from Lemma \ref{blow.up.base.change.lemm}, noting that the preimage of $S_0'$ under the flat map $\calG\x_SS'\to S'$ defines an effective Cartier divisor.
\xpf

\exam\label{neron.blow.exam}
Let $G_0\to \Spec(\bbZ)$ be a Chevalley group scheme, that is, a split reductive
group scheme with connected fibers \cite[\S6.4]{Co14} (note that some
authors call this a Demazure group scheme, keeping the term
Chevalley group scheme for split semisimple groups). Consider the base
change $G:=G_0\x_\bbZ\bbA^1_\bbZ$ to the affine line.
Let $S_0=\Spec(\bbZ)$ considered as the effective Cartier divisor defined
by the zero section of $S=\bbA^1_\bbZ$.
Let $P_0\subset G_0$ be a parabolic subgroup.
By Theorem \ref{neron.blow.theo} (3) and (5), the N\'eron blowup $\calG\to \bbA^1_\bbZ$ of $G$ in $P_0$ is a smooth, affine group scheme with connected fibers. 
In fact, it is an easy special case of the group schemes constructed in \cite[\S 4]{PZ13} and \cite[\S 2]{Lou}.
Let $\varpi$ denote a global coordinate on $\bbA^1_\bbZ$.
By Theorem \ref{neron.blow.theo} (6) the base changes have the following properties:
\begin{enumerate}
\item If $k$ is any field, then $\calG({k\pot{\varpi}})$ is the subgroup of those elements in $G(k\pot{\varpi})=G_0(k\pot{\varpi})$ whose reduction modulo $\varpi$ lies in $P_0(k)$.
\item If $p$ is any prime number, then $\calG_{\varpi\mapsto p}(\bbZ_p)$ is subgroup of those elements in $G_{\varpi\mapsto p}(\bbZ_p)=G_0(\bbZ_p)$ whose reduction modulo $p$ lies in $P_0(\bbF_p)$.
\end{enumerate}
In other words the respective base changes $\calG\x_{\bbA^1_\bbZ}\Spec(k\pot{\varpi})$ and $\calG\x_{\bbA^1_\bbZ,\varpi\mapsto p}\Spec(\bbZ_p)$ are parahoric group schemes in the sense of \cite{BT84}, cf.~also \cite[Cor.~4.2]{PZ13} and \cite[\S2.6]{Lou}.
Thus, the N\'eron blowup $\calG\to \bbA^1_\bbZ$ can be viewed as a family of parahoric group schemes.
\xexam
 
\subsection{Group structure on the exceptional divisor}
We continue with the notation of \S\ref{neron.blow.def.sec}.
In this subsection we take up the description of the exceptional
divisor from Proposition~\ref{exceptional divisor}, in the context
of group schemes: we explain the interplay between the ambient
group structure and the vector bundle structure on the exceptional
divisor.

In the sequel, for any scheme $X$ we write
$\Ga(X)=H^0(X,\calO_X)$ its ring of global functions.

\lemm \label{defi:inclusion in ker i times i}
Assume that $S$ is affine. Let $\calG$ be an $S$-group scheme
and $\calG_0:=\calG\x_SS_0$. Denote by $i\co \calG_0\hookto \calG$
the closed immersion and $\calK$ the corresponding ideal sheaf.
Let $m,\pr_1,\pr_2\co \calG\times_S\calG\to \calG$ be the multiplication
and the projections, with corresponding morphisms:
\[
\renewcommand{\arraystretch}{1.5}
\begin{array}{l}
m^\sharp,\pr_1^\sharp,\pr_2^\sharp
\co\Ga(\calG)\longto \Ga(\calG\times_S\calG) \\
(i\times i)^\sharp
\co\Ga(\calG\times_S\calG)\longto \Ga(\calG_0\times_{S_0} \calG_0).
\end{array}
\]
If $\delta:=m^\sharp-\pr_1^\sharp-\pr_2^\sharp$, we have
$\delta(H^0(\calG,\calK)) \subset\ker\left((i\times i)^\sharp\right)$.
\xlemm

\pf
Each of the maps $f\in\{m,\pr_1,\pr_2\}$ fits in a
commutative diagram:
\[
\begin{tikzcd}
\calG_0\times_{S_0}\calG_0 \ar[r,"i\times i"] \ar[d,"f"']
& \calG\times_S \calG \ar[d,"f"] \\
\calG_0 \ar[r,"i"] & \calG.
\end{tikzcd}
\]
Since $H^0(\calG,\calK)$ is the kernel of the map
$i^\sharp\co \Ga(\calG)\to \Ga(\calG_0)$, by
taking global sections we obtain
$((i\times i)^\sharp f^\sharp)(H^0(\calG,\calK))=0$.
\xpf

We recall that for a group scheme $G\to S$ with unit section
$e\co S\to G$, the Lie algebra $\Lie(G/S)$ is the $S$-group scheme
$\bbV(e^*\Omega^1_{G/S})$. 

\theo \label{theo:dilatations of subgroups}
With the notation of \S\ref{neron.blow.def.sec}, assume that
$G\to S$ is flat, locally finitely presented and $H\to S_0$ is flat,
regularly immersed in $G_0$. Let $\calG\to G$ be the dilatation
of $G$ in $H$ with exceptional divisor $\calG_0:=\calG\x_SS_0$.
Let $\calJ$ be the ideal sheaf of $G_0$ in $G$ and
$\calJ_H:=\calJ|_H$. Let $V$ be the restriction of the normal
bundle $\bbV(\calC_{H/G_0}\otimes \calJ_H^{-1})\to H$ along the
unit section $e_0\co S_0\to H$.
\begin{enumerate}
\item Locally over $S_0$, there is an exact sequence of
$S_0$-group schemes $1 \to V \to \calG_0 \to H\to 1$.
\item If $H$ lifts to a flat $S$-subgroup scheme of~$G$,
there is globally an exact, canonically split sequence
$1 \to V \to \calG_0 \to H\to 1$.
\item If $G\to S$ is smooth, separated and $\calG\to G$ is the
dilatation of the unit section of $G$, there is a canonical
isomorphism of smooth $S_0$-group schemes
$\calG_0 \isomto \Lie(G_0/S_0)\otimes\Norm_{S_0/S}^{-1}$
where $\Norm_{S_0/S}$ is the normal bundle of $S_0$ in $S$.
\end{enumerate}
\xtheo

\pf
(1) Let $\calF=\calC_{H/G_0}\otimes \calJ_H^{-1}$.
According to Proposition~\ref{exceptional divisor}(1), locally
over $S_0$ we have an isomorphism of $S_0$-schemes:
\[
\psi:\calG\times_G H \isomto \bbV(\calF).
\]
Let $K=\ker(\calG_0\to H)$. To obtain the exact
sequence of the statement, it is enough to prove that the
restriction of $\psi$ along the unit section
$e_0\co S_0\to H$ is an isomorphism of $S_0$-group schemes
\[
e_0^*\psi\co K\isomto V.
\]
For this we may localize further around a point of $S_0$, hence
assume that~$S$ and $S_0$ are affine and small enough so that
$\calJ$ is trivial.
Proving that $e_0^*\psi$ is a morphism of groups is equivalent to
checking an equality between two morphisms $K\times_{S_0} K\to V$.
Since $V=\Spec(\Sym(\calF_0))$ where $\calF_0:=e_0^*\calF$,
this is the same as checking an equality between two maps
of $\Ga(S_0)$-modules $H^0(S_0,\calF_0)\to \Ga(K\times_{S_0} K)$.
More precisely, since $K$ is affine we have
$\Ga(K\times_{S_0} K)=\Ga(K)\otimes_{\Ga(S_0)} \Ga(K)$ and what we have to
check is that $m^\sharp(x)=x\otimes 1+1\otimes x$ for all
$x\in H^0(S_0,\calF_0)$, with $m$ the multiplication
of $K$. That is, we want to prove that
\[
\delta(H^0(S_0,\calF_0))=0
\]
where $\delta:\Ga(K)\to \Ga(K\times_{S_0} K)$ is defined by
$\delta=m^\sharp-\pr_1^\sharp-\pr_2^\sharp$.

In order to prove this, let $\calI$ be the ideal sheaf of the
closed immersion $H\hookto G$, and let $f^*\calI$ be its preimage
as a module under the dilatation morphism $f\co \calG\to G$.
Consider the closed immersions
$K\hookto \calG_0$ and $i\co \calG_0\hookto \calG$, and the diagram:
\[
\begin{tikzcd}[row sep=15]
H^0(\calG,f^*\calI) \ar[r] \ar[dd] & \Ga(\calG) \ar[r,"\delta"] \ar[dd]
& \Ga(\calG\times_S\calG) \ar[d,"(i\times i)^\sharp"] \\
& & \Ga(\calG_0\times_{S_0} \calG_0) \ar[d]
\\
H^0(S_0,\calF_0) \ar[r] & \Ga(K) \ar[r,"\delta"] & \Ga(K\times_{S_0} K).
\end{tikzcd}
\]
We claim that the vertical map $H^0(\calG,f^*\calI)\to H^0(S_0,\calF_0)$
is surjective. To prove this, let $e_{\calG}\co S\to\calG$ be the
unit section of $\calG$ and $j\co S_0\hookto S$ be the closed immersion,
and decompose the said map as follows:
\[
H^0(\calG,f^*\calI) \stackrel{e_{\calG}^*}{\longto}
H^0(S,e^*\calI) \stackrel{j^*}{\longto} H^0(S_0,j^*e^*\calI)
\longto H^0(S_0,\calF_0).
\]
The first map is surjective because it has the section
$\sigma^*$ where $\sigma\co \calG\to S$ is the structure
map. The second map is surjective because it is obtained by
taking global sections on the affine scheme $S$ of the surjective
map of sheaves $e^*\calI\to j_*j^*e^*\calI$. To show that the third
map is surjective, start from the surjection of sheaves
$\calI|_H\to\calF$. Since pullback is right exact, this gives
rise to a surjection $j^*e^*\calI=e_0^*\calI|_H\to e_0^*\calF=\calF_0$.
Taking global sections on the affine scheme $S_0$ we obtain
the desired surjection.

Now let $\calI\calO_{\calG}=f^{-1}\calI\cdot\calO_{\calG}$
be the preimage of $\calI$ as an ideal.
Note that by property of the dilatation, we have
$\calI\calO_{\calG}=\calJ\calO_{\calG}=:\calK$. Therefore, according to
Lemma~\ref{defi:inclusion in ker i times i}, we have
\[
\delta(H^0(\calG,\calI\calO_{\calG}))=
\delta(H^0(\calG,\calK))\subset\ker ((i\times i)^\sharp).
\]
Precomposing with the surjection $f^*\calI\to\calI\calO_{\calG}$,
we find that $H^0(\calG,f^*\calI)$ is mapped into
$\ker(i\times i)^\sharp$ by $\delta$.
As a result $H^0(\calG,f^*\calI)$ goes to zero
in $\Ga(K\times_{S_0} K)$. Since
$H^0(\calG,f^*\calI)\to H^0(S_0,\calF_0)$ is surjective,
the commutativity of the diagram implies that
$\delta(H^0(S_0,\calF_0))=0$ in $\Ga(K\times_{S_0} K)$, as desired.
Hence $e_0^*\psi:K\to V$ is an isomorphism of groups. This proves (1).

\smallskip

\noindent (2)
If $\tilde H\subset G$ is a flat $S$-subgroup scheme lifting $H$,
we have a transversal intersection $H=\tilde H\cap G_0$.
By Proposition~\ref{exceptional divisor}(3), the preceding
construction of the short exact sequence can be performed globally
over $S_0$. Moreover, by
the universal property of the dilatation the map
$\tilde H\to G$ lifts to a map $\tilde H\to \calG$. In restriction
to $S_0$ this splits the short exact sequence previously obtained.

\smallskip

\noindent (3)
Finally if $G\to S$ is smooth, the unit section is a regular
immersion with conormal sheaf $\omega_{G/S}=e^*\Omega^1_{G/S}$.
In restriction to $S_0$ the group $V$ is the Lie algebra whence
the canonical isomorphism
$\calG_0\isomto \Lie(G_0/S_0)\otimes\Norm_{S_0/S}^{-1}$.
\xpf

\rema
In the situation of Theorem \ref{theo:dilatations of subgroups} (2),
the group $H$ acts by conjugation on
$V=\bbV(e_0^*\calC_{H/G_0}\otimes \calJ_{S_0}^{-1})$.
We checked on examples that this additive action is linear, and is
in fact none other than the ``adjoint'' representation of $H$ on
its normal bundle as in \cite[Exp.~I, Prop.~6.8.6]{SGA3.1}. 
To recall what this representation is, note that the normal sheaf
$\calJ_{S_0}^{-1}$ comes from the base
and is endowed with the trivial action; for simplicity we describe
the action Zariski locally on $S_0$ and omit it from the notation.
We start from the conormal
sequence of the inclusions $\{1\}\subset H\subset G_0$:
\[
\begin{tikzcd}[column sep=15]
0 \ar[r,dashed] & e_0^*\calC_{H/G_0} \ar[r] & \omega_{G_0} \ar[r]
& \omega_H \ar[r] & 0.
\end{tikzcd}
\]
The sequence is exact on the left if $\{1\}\to H$ is a regular
immersion. Taking the associated vector bundles, we obtain an
exact sequence of $S_0$-group schemes:
\[
\begin{tikzcd}[column sep=15]
0 \ar[r] & \Lie(H/S_0) \ar[r] & \Lie(G_0/S_0) \ar[r] &
\bbV(e_0^*\calC_{H/G_0}) \ar[r,dashed] & 0.
\end{tikzcd}
\]
The middle term $\Lie(G_0/S_0)$ supports the adjoint action of $G_0$.
The restricted action of $H$ leaves stable the terms $\Lie(H/S_0)$
and $\bbV(e_0^*\calC_{H/G_0})$, cf.~\cite[Exp.~III, Lem.~4.25]{SGA3.1}. The former is the adjoint action of $H$
on its Lie algebra, and the latter is the action on the normal
bundle along the unit section.
\xrema

\subsection{N\'eron blowups as syntomic sheaves}
\label{syntomic.sheaf.sec}
We continue with the notation of \S\ref{neron.blow.def.sec}. 
Additionally assume that $j\co S_0\hookto S$ is an effective Cartier divisor,
that $G\to S$ is a flat, locally finitely presented group scheme
and that $H\subset G_0:=G\times_SS_0$ is a flat, locally finitely
presented closed
$S_0$-subgroup scheme. In this context, there is another
viewpoint on the dilatation $\calG$ of $G$ in $H$, namely as the
kernel of a certain map of syntomic sheaves.

To explain this, let $f\co G_0\to G_0/H$ be the morphism to
the fppf quotient sheaf, which by Artin's theorem
(\cite[Cor.~6.3]{Ar74}) is representable by an algebraic
space. By the structure theorem for algebraic group schemes
(see \cite[Exp.~VII$_{\on{B}}$, Cor.~5.5.1]{SGA3.1})
the morphisms $G\to S$ and $H\to S_0$ are syntomic. Since
$f\co G_0\to G_0/H$ makes $G_0$ an $H$-torsor, it follows
that $f$ is syntomic also.

\lemm \label{lemma:short exact sequence}
Let $S_{\syn}$ be the small syntomic site of $S$. Let
$\eta\co G\to j_*j^*G$ be the adjunction map in the category
of sheaves on $S_{\syn}$ and consider the composition
$v=(j_*f)\circ \eta$:
\[
\begin{tikzcd}
G \ar[r,"\eta"] & j_*j^*G=j_*G_0 \ar[r,"j_*f"] & j_*(G_0/H).
\end{tikzcd}
\]
Then the dilatation $\calG\to G$ is the kernel of $v$.
More precisely, we have an exact sequence of sheaves of
pointed sets in $S_{\syn}$:
\[
\begin{tikzcd}[column sep=15]
1 \ar[r] & \calG \ar[r] & G \ar[r,"v"] & j_*(G_0/H) \ar[r] & 1.
\end{tikzcd}
\]
If $G\to S$ and $H\to S_0$ are smooth, then the sequence is
exact as a sequence of sheaves on the small \'etale site of $S$.
\xlemm

\pf
That $\calG\to G$ is the kernel of $v$ follows directly from the universal
property of the dilatation, restricted to syntomic $S$-schemes.
It remains to prove that the map of sheaves $v$ is surjective.
It is enough to prove that both maps~$\eta$ and $j_*f$ are
surjective. For $\eta$ this is because if $T\to S$ is a syntomic
morphism, any point $t\co T\to j_*G_0$ lifts tautologically to $G$
after the syntomic refinement $T'=G\times_ST\to T$. For $j_*f$,
we start from a syntomic morphism $T\to S$ and a point
$t\co T\to j_*(G_0/H)$, that is a morphism $T_0\to G_0/H$.
Using that $G_0\to G_0/H$ is syntomic, we can find as before
a syntomic refinement $T'_0\to T_0$ and a lift $T'_0\to G_0$.
Using that syntomic coverings lift across closed immersions
(see
\cite[\href{https://stacks.math.columbia.edu/tag/04E4}{04E4}]{StaPro}),
there is a syntomic covering $T''\to T$ such that $T''_0$
refines $T'_0$. This provides a lift of $t$ to $j_*G_0$.

Finally if $G\to S$ and $H\to S_0$ are smooth, the existence
of \'etale sections for smooth morphisms (\cite[Cor.~17.16.3]{EGA4.4})
and the possibility to
lift \'etale coverings across closed immersions
(see \cite[\href{https://stacks.math.columbia.edu/tag/04E4}{04E4}]{StaPro}
again) show that the sequence is exact also in the \'etale topology.
\xpf

\section{Applications}\label{applications.sec}
Here we give two applications in cohomological degree $0$ and $1$ of the theory developed so far: integral points and torsors. 
In \S\ref{Moy.Prasad.sec} we consider integral points of N\'eron blowups and discuss the isomorphism relating the graded pieces of the congruent filtration of $G$ to the graded pieces of its Lie algebra $\frakg$.
In \S\ref{Torsors.blowup.sec} we discuss torsors under N\'eron blowups and apply this in \S\ref{bundles.curves.sec} to the construction of level structures on moduli stacks of $G$-bundles on curves, and in \S\ref{applications.integral.models.sec} to the construction of integral models of moduli stacks of shtukas.

\subsection{Integral points and the Moy-Prasad
isomorphism}\label{Moy.Prasad.sec}
In this subsection we prove an isomorphism
describing the graded pieces of the filtration by congruence subgroups
on the integral points of reductive group schemes. For the benefit of the interested reader,
we provide comments on the literature on this topic in
Remark~\ref{literature.on.MP.remark} below.

We start with the following lemma.

\lemm \label{lemma:congruence}
Let $\calO$ be a ring and $\pi\subset \calO$ an invertible ideal
such that $(\calO,\pi)$ is a henselian pair. Let $G$ be a smooth,
separated $\calO$-group scheme and $\calG\to G$ the
dilatation of the trivial subgroup over $\calO/\pi$. If either
$\calO$ is local or~$G$ is affine, then the exact sequence of
Lemma~\ref{lemma:short exact sequence}
induces an exact sequence of groups:
\[
1\longto \calG(\calO) \longto G(\calO) \longto G(\calO/\pi) \longto 1.
\]
\xlemm

\pf
Write $S=\Spec(\calO)$, $S_0=\Spec(\calO/\pi)$ and set $G_0=G\x_SS_0$,
$\calG_0=\calG\x_SS_0$. Consider the short exact sequence of
Lemma~\ref{lemma:short exact sequence} on the \'etale site and take the
global sections over $S$. It is then enough to prove that the map
$G(\calO) \to G(\calO/\pi)$ is surjective. For this start with an
$(\calO/\pi)$-point of $G$, i.e. a section $u_0:S_0\to G_0$ to the
map $G_0\to S_0$. If either $\calO$ is local or~$G$ is affine,
$u_0$ factors through an open affine subscheme $U\subset G$.
In this situation, the classical existence result for lifting of
sections for smooth schemes over a henselian local ring
(as in for example \cite[2.3/5]{BLR90})
extends to henselian pairs, see \cite[Thm.~1.8]{Gr72}.
In this way we see that $u_0$ lifts to a section $u$ of $G\to S$.
\xpf

\rema \rm
The same proof gives a similar result with the dilatation
of a split smooth unipotent closed subgroup $H\subset G_0$,
with the group $G(\calO/\pi)$ replaced by the pointed set
$(G_0/H)(\calO/\pi)$. Indeed, such a group $H$ is obtained
by successive extensions of the additive group $\bbG_{a,S_0}$
over $S_0=\Spec(\calO/\pi)$. Since $S_0$ is affine, the
(\'etale or syntomic) cohomology of $\bbG_{a,S_0}$ vanishes, being
the coherent cohomology of $\calO_{S_0}$. Using induction one
concludes that $H^1(S_0,H)$ is trivial. Now starting from a section
$u_0$ of $G_0/H\to S_0$ and pulling back the map $G_0\to G_0/H$
along $u_0$, we obtain an $H$-torsor $G_0\times_{G_0/H}S_0\to S_0$.
By the previous remarks this torsor has a section $v_0$.
The latter lifts to a section of $G$ by the same argument as in the
proof of the lemma.
\xrema

\theo \label{theoisomp}
Let $r,s$ be integers such that $0\le r/2\le s\le r$.
Let $(\calO,\pi)$ be a henselian pair where $\pi\subset \calO$ is
an invertible ideal. Let $G$ be a smooth, separated
$\calO$-group scheme. Let $G_r$ the $r$-th iterated dilatation of
the unit section and $\frakg_r$ its Lie algebra. If $\calO$ is
local or $G$ is affine, there is a canonical  isomorphism:
\[
G_s(\calO)/G_r(\calO) \isomto \frakg_s(\calO)/\frakg_r(\calO).
\]
\xtheo

\pf
The $(r-s)$-th iterated dilatation of $G_s$ is naturally $G_r$.
But as we observed in
Proposition~\ref{prop:dilatation of thick subscheme},
the group scheme $G_r$ can also be seen as the dilatation
of $\{1\}$ in $(\calO,\pi^r)$.
For an integer $n\geq 0$, we write $\calO_n:=\calO/\pi^n$.
Putting these remarks together,
the previous lemma applied to the dilatation of the group scheme
$G=G_s$ with respect to $(\calO,\pi^{r-s})$ yields an isomorphism:
\begin{equation}\label{quotient.integral.pts.congruent}
G_s(\calO)/G_r(\calO) \isomto G_s(\calO_{r-s}).
\end{equation}
We now consider the statement of the theorem. If $s=0$
we have $r=0$, hence left-hand side and right-hand side
are equal to $\{1\}$ and the result is clear. Therefore we may
assume that $s>0$.
Theorem~\ref{theo:dilatations of subgroups}
applied to the dilatation of~$\{1\}$ in $(\calO,\pi^s)$
provides a canonical isomorphism
\[
{G_s}|_{\calO_s} \isomto \Lie(G|_{\calO_s})\otimes\Norm_{\calO_s/\calO}^{-1}.
\]
Since the $\Lie$ algebra of a vector bundle $\bbV(\scrE)$ is
canonically isomorphic to $\bbV(\scrE)$ itself
(\cite[Exp.~II, Ex.~4.4.2]{SGA3.1}), taking Lie algebras on
both sides we deduce a canonical isomorphism
\[
{G_s}|_{\calO_s} \isomto {\frakg_s}|_{\calO_s}.
\]
Since $r-s\le s$, the ring $\calO_{r-s}$ is an $\calO_s$-algebra
and we can take $\calO_{r-s}$-valued points in the previous
isomorphism to obtain:
\[
G_s(\calO_{r-s}) \isomto \frakg_s(\calO_{r-s}).
\]
Using \eqref{quotient.integral.pts.congruent} once for $G_s$ and once
for $\frakg_s$, we end up with
\[
G_s(\calO)/G_r(\calO) \isomto \frakg_s(\calO)/\frakg_r(\calO),
\]
which is the desired canonical isomorphism.
\xpf

\rema \label{literature.on.MP.remark}
In the literature on integral points of reductive groups over
non-archimedean local fields, results similar to the isomorphism
of Theorem~\ref{theoisomp} appeared with restrictions on the
indices $r,s$, on the group schemes involved or on the ground
ring. See for instance
\cite[Prop.~6 (b)]{Ser68}, \cite[Prop.~3.9 and 3.10]{Ne99}
for the multiplicative group, \cite[p.~442 line 1]{Ho77}, \cite{Mo82}, \cite[p.~22]{BK93} for general linear
groups, \cite[p.~337]{Sec04} for general linear groups over division
algebras, and \cite[\S2]{PR84}, \cite[\S2]{MP94}, \cite{MP96},
\cite[\S1]{Ad98}, \cite[\S1]{Yu01} for general reductive groups.
In these examples, the isomorphisms are defined at the level of
integral points using ad hoc explicit formulas. These isomorphisms
are sometimes called Moy-Prasad isomorphisms as a tribute to
\cite{MP94}, read \cite[0.4]{Yu15} and \cite[p.~242 lines 18-19]{De02}
for informations. In the case of an affine, smooth group scheme
over a discrete valuation ring, the isomorphism of
Theorem~\ref{theoisomp} appears without proof in 
\cite[proof of Lemma 2.8]{Yu15}.
\xrema

\subsection{Torsors and level structures}\label{Torsors.blowup.sec}
In this subsection we adopt notations more specific to the study
of torsors over curves. 
Let $X$ be a scheme, and let $N\subset X$ be an effective
Cartier divisor.
Let $G\to X$ be a smooth, finitely presented group scheme, and let $H\subset G|_N$ be an $N$-smooth closed subgroup. 
We denote by $\calG\to G$ the N\'eron blowup of $G$ in $H$ (over $N$) which is a smooth, finitely presented $X$-group scheme by Theorem \ref{neron.blow.theo}.

For a scheme $T\to X$, let $BG(T)$ (resp.~$B\calG(T)$) denote the groupoid of right $G$-torsors on $T$ in the fppf topology. 
Here we note that every such torsor is representable by a smooth algebraic space (of finite presentation), and hence admits sections \'etale locally.
Whenever convenient we may therefore work in the \'etale topology as opposed to the fppf topology.

Pushforward of torsors along $\calG\to G$ induces a morphism of contravariant functors $\Sch_X\to \Gpds$ given by
\begin{equation}\label{Neron.blow.stacks.eq}
B\calG\to BG, \;\; \calE\mapsto \calE\x^\calG G.
\end{equation}

\defi\label{neron.blow.groupoid.defi}
For a scheme $T\to X$, let $B(G,H,N)(T)$ be the groupoid whose objects are pairs $(\calE,\be)$ where $\calE\to T$ is a right fppf $G$-torsor and $\be$ is a section of the fppf quotient 
\[
\big(\calE|_{T_N}\big/H|_{T_N}\big)\to T_N,
\]
where $T_N:=T\x_XN$, i.e., $\be$ is a reduction of $\calE|_{T_N}$ to an $H$-torsor.
Morphisms $(\calE,\be)\to (\calE',\be')$ are given by isomorphisms of torsors $\varphi\co \calE\cong \calE'$ such that $\bar \varphi\circ\be=\be'$ where $\bar \varphi$ denotes the induced map on the quotients.
Note that if $T_N=\varnothing$, then there is no condition on the compatibility of $\be$ and $\be'$.
\xdefi

Each of the contravariant functors $\Sch_X \to \Gpds$ induced by $B\calG$, $BG$ and $B(G,H,N)$ defines a stack over $X$ in the fppf topology. 
We call $B(G,H,N)$ the {\it stack of $G$-torsors equipped with level-$(H,N)$-structures}. 

\lemm
The map \eqref{Neron.blow.stacks.eq} factors as a map of $X$-stacks
\begin{equation}\label{neron.blow.factor.eq}
B\calG\to B(G,H,N)\to BG,
\end{equation}
where the second arrow denotes the forgetful map.
\xlemm
\pf
By Lemma \ref{Neron.blow.lemm} (1) the map $\calG|_N\to G|_N$ factors as $\calG|_N\to H\subset G|_N$. 
Thus, given a $\calG$-torsor $\calE\to T$ we get the $H$-equivariant map 
\[
\calE\x^{\calG|_{T_N}}H|_{T_N} \; \subset\; \calE\x^{\calG|_{T_N}}G|_{T_N}.
\]
Passing to the fppf quotient for the right $H$-action defines the section $\be_{\on{can}}$.
The association $\calE\mapsto (\calE\x^\calG G, \be_{\on{can}})$ induces the desired map $B\calG\to B(G,H,N)$.
\xpf

\prop \label{Neron.blow.equivalence.cor}
The map \eqref{neron.blow.factor.eq} induces an equivalence of contravariant functors $\Sch_X^{N\text{-}\reg}\to \Gpds$ given by
\[
B\calG\overset{\simeq}{\longto} B(G,H,N), \;\; \calE\mapsto (\calE\x^\calG G,\, \be_{\on{can}}).
\]
\xprop
\pf 
For $T\to X$ in $\Sch_X^{N\text{-}\reg}$, we need to show that $B\calG(T)\to B(G,H,N)(T)$ is an equivalence of groupoids. 
Since $\calG\to X$ is smooth, in particular flat, its formation commutes with base change along $T\to X$ by Theorem \ref{neron.blow.theo}.
Hence, we may reduce to the case where $T=X$.
Now recall from Lemma \ref{lemma:short exact sequence} the exact sequence of sheaves of pointed sets on the \'etale site of $X$,
\[
\begin{tikzcd}[column sep=15]
1 \ar[r] & \calG \ar[r] & G \ar[r] & j_*(G|_N/H) \ar[r] & 1,
\end{tikzcd}
\]
where $j\co N\subset X$ denotes the inclusion. 
The desired equivalence is a consequence of \cite[Chap.~III, \S 3.2, Prop.~3.2.1]{Gi71}. 
A quasi-inverse is given by pulling back a section $\be\co X\to  j_*(G|_N/H)$ along the $\calG$-torsor $G \to j_*(G|_N/H)$. 
Here we have used that by smoothness of the group schemes $G\to X$, $\calG\to X$, $H\to N$ and consequently of the quotient $G|_N/H$ we are allowed to work with the \'etale topology as opposed to the fppf topology. 
\xpf

\subsubsection{Level structures on moduli stacks of bundles on curves}\label{bundles.curves.sec}
We continue with the notation, and additionally assume that $X$ is a smooth, projective, geometrically irreducible curve over a field $k$, and that $G\to X$ and hence $\calG\to X$ is affine. 

Let $\Bun_{G}:=\Res_{X/k}BG$ (resp.~$\Bun_{\calG}:=\Res_{X/k}B\calG$)
be the moduli stack of $G$-torsors (resp.~$\calG$-torsors) on $X$;
here $\Res_{X/k}$ stands for the Weil restriction along $X\to \Spec(k)$.
This is a quasi-separated, smooth algebraic stack locally of finite type
over $k$, cf.~e.g.~\cite[Prop.~1]{He10} or \cite[Thm.~2.5]{AH19}. Similarly, let
$\Bun_{(G,H,N)}:=\Res_{X/k}B(G,H,N)$ be the stack parametrizing $G$-torsors
over $X$ with level-$(H,N)$-structures as in
Definition \ref{neron.blow.groupoid.defi}.

\theo\label{applications.bun.level.theo}
The map \eqref{neron.blow.factor.eq} induces equivalences of contravariant functors $\Sch_k\to \Gpds$ given by
\[
\Bun_\calG\overset{\cong}{\longto} \Bun_{(G,H,N)}, \;\; \calE\mapsto (\calE\x^\calG G,\, \be_{\on{can}})
\]
\xtheo
\pf 
For any $k$-scheme $T$, the projection $X\x_kT \to X$ is flat, and hence defines an object in $\Sch_X^{N\on{-}\reg}$.
The theorem follows from Proposition \ref{Neron.blow.equivalence.cor}.
\xpf

\exam \label{applications.exam.LS}
If $H=\{1\}$ is trivial, then $\Bun_{(G,H,N)}$ is the moduli stack of $G$-torsors on $X$ with level-$N$-structures.
If $G\to X$ is split reductive, if $N$ is reduced and if $H$ a parabolic subgroup in $G|_N$, then $\Bun_{(G,H,N)}$ is the moduli stack of $G$-torsors with quasi-parabolic structures in the sense of Laszlo-Sorger \cite{LS97}, cf.~\cite[\S 2.a.]{PR10} and \cite[\S1, Exam.~(2)]{He10}.
\xexam

We end this subsection by discussing Weil uniformizations. 
Let $|X|\subset X$ be the set of closed points, and let $\eta \in X$ be the generic point.
We denote by $F=\kappa(\eta)$ the function field of $X$.
For each $x\in |X|$, we let $\calO_x$ be the completed local ring at $x$ with fraction field $F_x$ and residue field $\kappa(x)=\calO_x/\frakm_x$. 
Let $\bbA:=\bigsqcap'_{x\in |X|}F_x$ be the ring of adeles with subring of integral elements $\bbO=\bigsqcap_{x\in |X|}\calO_x$.
As in \cite[Lem.~1.1]{N06} or \cite[Rem.~8.21]{Laf18} one has the following result. 

\prop\label{applications.Weil.uniformization.prop}
Assume that $k$ is either a finite field or a separably closed field,
and that $G\to X$ has connected fibers. 
Then there is an equivalence of groupoids
\begin{equation}\label{applications.Weil.uniformization.eq}
\Bun_G(k)\;\simeq\; \bigsqcup_{\ga}  G_\ga(F)\big\bsl\big( G_\ga(\bbA)/G(\bbO)\big),
\end{equation}
where $\ga$ ranges over $\ker^1(F,G):=\ker\big(H^1_{\et}(F,G)\to \bigsqcap_{x\in |X|}H^1_{\et}(F_x,G)\big)$, and where $G_\ga$ denotes
the associated pure inner form of $G|_F$.
The identification \eqref{applications.Weil.uniformization.eq} is functorial in $G$ among maps of $X$-group schemes which are isomorphisms in the generic fibre. 
\xprop
\pf 
Under our assumptions, Lang's lemma implies that $H^1_\et(\calO_x,G)$ is trivial for all $x\in |X|$: 
use that $H^1_\et(\kappa(x),G)$ is trivial because $G|_{\kappa(x)}$ smooth, affine, connected and $\kappa(x)$ is either finite or separably closed; then an approximation argument as in e.g.~\cite[Lem.~A.4.3]{RS20}. 
In particular, for every $G$-torsor $\calE\to X$ the class of its generic fibre $[\calE|_F]$ lies in $\ker^1(F,G)$. 
For each $\ga\in \ker^1(F,G)$, we fix a $G$-torsor $\calE_\ga^0\to \Spec(F)$ of class $\ga$. 
We denote by $G_\ga$ its group of automorphisms which is an inner form of $G$.
We also fix an identification $G_\ga(F_x)=G(F_x)$ for all $x\in |X|$, $\ga\in  \ker^1(F,G)$. 
In particular, $G_\ga(\bbA)=G(\bbA)$ so that the right hand quotient in \eqref{applications.Weil.uniformization.eq} is well-defined.
Now consider the groupoid
\[
\Sig_\ga:=\big\{(\calE,\delta,(\epsilon_{x})_{x\in |X|})\;\big|\; \delta\co \calE|_F\simeq \calE_\ga^0,\;\; \epsilon_x\co \calE^0\simeq \calE|_{\calO_x}  \big\}.
\] 
For each $x\in |X|$, we have
\[
g_x:=\delta|_{F_x}\circ \epsilon_x|_{F_x}\in \Aut(\calE_\ga^0|_{F_x})=G_\ga(F_x)=G(F_x),
\]
and further $g_x\in G(\calO_x)$ for almost all $x\in |X|$. 
Thus, the collection $(g_x)_{x\in |X|}$ defines a point in $G(\bbA)=G_\ga(\bbA)$. 
In this way, we obtain an $G_\ga(F)\x G(\bbO)$-equivariant map $\pi_\ga\co \Sig_\ga\to G_\ga(\bbA)$, and thus a commutative diagram of groupoids
\[
\begin{tikzpicture}[baseline=(current  bounding  box.center)]
\matrix(a)[matrix of math nodes, 
row sep=1.5em, column sep=2em, 
text height=1.5ex, text depth=0.45ex] 
{ 
\bigsqcup_\ga\Sig_\ga& \bigsqcup_\ga G_\ga(\bbA) \\ 
\Bun_G& \bigsqcup_\ga G_\ga(F)\big\bsl\big(G_\ga(\bbA)/G(\bbO )\big).\\}; 
\path[->](a-1-1) edge node[above] {$\sqcup_\ga\pi_\ga$} (a-1-2);
\path[->](a-2-1) edge[dashed] node[above] {} (a-2-2);
\path[->](a-1-1) edge node[right] {} (a-2-1);
\path[->](a-1-2) edge node[right] {} (a-2-2);
\end{tikzpicture}
\] 
As the vertical maps are disjoint unions of $G_\ga(F)\x G(\bbO)$-torsors, the dashed arrow is fully faithful. 
Hence, it suffices to show that it is a bijection on isomorphism classes, i.e., a bijection of sets.
We construct an inverse of the dashed arrow as follows: Given a representative $(g_x)_{x\in |X|}\in G_\ga(\bbA)=G(\bbA)$ of some class, there is a non-empty open subset $U\subset X$ such that $g_x\in G(\calO_x)$ for all $x\in |U|$, and such that $\calE_\ga^0$ is defined over $U$. 
Let $X\bsl U=\{x_1,\ldots, x_n\}$ for some $n\geq 0$. We define the associated $G$-torsor by gluing the torsor $\calE^0_\ga$ on $U$ with the trivial $G$-torsor on  
\[
\Spec(\calO_{x_1})\sqcup\ldots\sqcup \Spec(\calO_{x_n})
\]
using the elements $g_{x_1},\ldots, g_{x_n}$ and the identification $G_\ga(F_x)=G(F_x)$. 
The gluing is justified by the Beauville-Laszlo lemma \cite{BL95}, or alternatively \cite[Lem.~5]{He10}. 
This shows \eqref{applications.Weil.uniformization.eq}. 
From the construction of the map $\sqcup_\ga\pi_\ga$, one sees that \eqref{applications.Weil.uniformization.eq} is functorial in $G$ among generic isomorphisms. 
\xpf

Note that $N$ defines an effective Cartier divisor on $\Spec(\bbO)$ so that the map of groups $\calG(\bbO)\to G(\bbO)$ is injective. 
As subgroups of $G(\bbO)$ we have
\begin{equation}\label{application.subgroup.eq}
\calG(\bbO)=\ker\big(G(\bbO)\to G(\calO_N)\to G(\calO_N)/H(\calO_N)\big),
\end{equation}
where $\calO_N$ denotes the ring of functions on $N$ viewed as a quotient ring $\bbO_X\to \calO_N$.

\coro \label{applications.Weil.uniformization.coro}
Under the assumptions of Proposition \ref{applications.Weil.uniformization.prop}, the N\'eron blowup $\calG\to X$ is smooth, affine with connected fibers by Theorem \ref{neron.blow.theo}, and there is a commutative diagram of groupoids
\[
\begin{tikzpicture}[baseline=(current  bounding  box.center)]
\matrix(a)[matrix of math nodes, 
row sep=1.5em, column sep=2em, 
text height=1.5ex, text depth=0.45ex] 
{ 
\Bun_\calG(k)& \bigsqcup_\ga G_\ga(F)\big\bsl\big(G_\ga(\bbA)/\calG(\bbO )\big) \\ 
\Bun_G(k)& \bigsqcup_\ga G_\ga(F)\big\bsl\big(G_\ga(\bbA)/G(\bbO )\big),\\}; 
\path[->](a-1-1) edge node[above] {$\simeq$} (a-1-2);
\path[->](a-2-1) edge node[above] {$\simeq$} (a-2-2);
\path[->](a-1-1) edge node[right] {} (a-2-1);
\path[->](a-1-2) edge node[right] {} (a-2-2);
\end{tikzpicture}
\] 
identifying the vertical maps as the level maps.
\xcoro

\rema
Let $G|_F$ be reductive.
\begin{enumerate}
\item If $k$ is algebraically closed, then $F$ is $C_1$ by Tsen's theorem and in particular of cohomological dimension $\leq 1$, cf.~\cite[II.3]{Ser65}. 
In this case, $H^1(F,G)$, and hence $\ker^1(F,G)$, is trivial by \cite[8.6]{BS68}.

\item If $k$ is a finite field, then $\ker^1(F,G)$ is dual to $\ker^1(F,Z(\hat G)(\bar \bbQ))$ where $Z(\hat G)$ denotes the center of the Langlands dual group $\hat G$, formed, say, over $\bar \bbQ$, cf.~\cite{Ko84, Ko86} and  \cite{NQT11} for global fields of positive characteristic.
In particular, if $G|_F$ is either simply connected or split reductive, then $\ker^1(F,G)$ is trivial, cf.~also \cite[Rem.~8.21, 12.2]{Laf18}.
\end{enumerate}
\xrema

\subsubsection{Integral models of moduli stacks of shtukas}\label{applications.integral.models.sec}
Here we point out that Theorem \ref{applications.bun.level.theo} immediately applies to construct certain integral models of moduli stacks of shtukas. 
We proceed with the notation of \S\ref{bundles.curves.sec}, and additionally assume that $k$ is a finite field.
Our presentation follows \cite[\S\S1-2]{Laf18}.

For any partition $I=I_1\sqcup \ldots\sqcup I_r$, $r\in \bbZ_{\geq 1}$ of a finite index set, the {\em moduli stack of iterated $G$-shtukas} is the contravariant functor of groupoids $\Sch_k\to \Gpds$ given by
\begin{equation}\label{shtuka.def}
\Sht_{G,I_\bullet}\defined \Big\{ \calE_r\overset{\al_r}{\underset{I_r}{\dashrightarrow}}\calE_{r-1}\overset{\al_r}{\underset{I_{r-1}}{\dashrightarrow}}\ldots\overset{\al_2}{\underset{I_2}{\dashrightarrow}} \calE_1\overset{\al_1}{\underset{I_1}{\dashrightarrow}} \calE_0={^\tau\calE_r} \Big\},
\end{equation}
where ${^\tau\calE}:= (\id_X\x \Frob_{T/k})^*\calE$ denotes the pullback under the relative Frobenius $\Frob_{T/k}$. 
Here the dashed arrows in \eqref{shtuka.def} indicate that the maps $\al_j$ between $G$-bundles are rationally defined.
More precisely, $\Sht_{G,I_\bullet}(T)$ classifies data $\big((\calE_j)_{j=1,\ldots, r}, \{x_i\}_{i\in I}, (\al_j)_{j=1,\ldots, r}\big)$
where $\calE_j\in \Bun_{G}(T)$ are torsors, $\{x_i\}_{i\in I}\in X^I(T)$ are points, and
\[
\al_j\co \calE_j|_{X_T\bsl (\cup_{i\in I_j}\Ga_{x_i})}\to \calE_{j-1}|_{X_T\bsl (\cup_{i\in I_j}\Ga_{x_i})}
\]
are isomorphisms of torsors. 
Here $\Ga_{x_i}\subset X_T$ denotes the graph of $x_i$ viewed as a relative effective Cartier divisor on $X_T\to T$. 
We have a forgetful map $\Sht_{G,I_\bullet}\to X^I$. 
Similarly, we have the moduli stack $\Sht_{\calG,I_\bullet}\to X^I$ defined by replacing $G$ with $\calG$.
By \cite{Var04} for split reductive groups and by \cite[Thm.~3.15]{AH19} for general smooth, affine group schemes both stacks are ind-(Deligne-Mumford) stacks which are ind-(separated and of locally finite type) over $k$.
Furthermore, pushforward of torsors along $\calG\to G$ induces a map of $X^I$-stacks
\begin{equation}\label{integral.model.eq2}
\Sht_{\calG,I_\bullet}\to \Sht_{G,I_\bullet},
\end{equation}
cf.~\cite{Br}. We also consider the {\em moduli stack of iterated $G$-shtukas with level-$(H,N)$-structures},
\[
\Sht_{(G,H,N),I_\bullet}\to X^I,
\] 
i.e., $\Sht_{(G,H,N),I_\bullet}(T)$ classifies data 
\[
\big((\calE_j,\be_j)_{j=1,\ldots, r}, \{x_i\}_{i\in I}, (\al_j)_{j=1,\ldots, r}\big),
\]
where $(\calE_j,\be_j)\in \Bun_{(G,H,N)}(T)$ are $G$-torsors with a level-$(H,N)$-structure, $\{x_i\}_{i\in I}\in X^I(T)$ are points, and
\begin{equation}\label{integral.model.eq}
\al_j\co (\calE_j,\be_j)|_{X_T\bsl (\cup_{i\in I_j}\Ga_{x_i})}\to (\calE_{j-1},\be_{j-1})|_{X_T\bsl (\cup_{i\in I_j}\Ga_{x_i})}
\end{equation}
are maps of $G$-torsors with a level-$(H,N)$-structure where $(\calE_0,\be_0):=({^\tau\calE}_{r},{^\tau\be}_{r})$. 
We have a forgetful map of $X^I$-stacks
\begin{equation}\label{integral.model.eq3}
\Sht_{(G,H,N),I_\bullet}\to \Sht_{G,I_\bullet}.
\end{equation}

\coro \label{applications.shtuka.model.theo}
Let $G, H, N$, $\calG$ and $I=I_1\sqcup\ldots\sqcup I_r$ be as above. Then the equivalence in Theorem~\ref{applications.bun.level.theo} induces an equivalence of $X^I$-stacks
\[
\Sht_{\calG,I_\bullet}\;\overset{\cong}{\longto}\;\Sht_{(G,H,N),I_\bullet},
\]
which is compatible with the maps \eqref{integral.model.eq2} and \eqref{integral.model.eq3}.
\xcoro

Loosely speaking, $\Sht_{\calG,I_\bullet}\cong \Sht_{(G,H,N),I_\bullet}$ over $X^I$ is an integral model for $\Sht_{(G,H,N),I_\bullet}|_{(X\bsl N)^I}$, which however needs modification outside the case of parahoric group schemes $\calG\to X$.
Concretely, if the characteristic places $x_i$, $i\in I$ of the shtuka divide the level $N$, then there is simply no compatibility condition on the $\be_j$'s in \eqref{integral.model.eq}. 
Consequently, the fibers of the map \eqref{integral.model.eq2}, resp.~\eqref{integral.model.eq3} over such places are by \cite[Thm.~3.20]{Br} certain quotients of positive loop groups. 
In particular, these fibers are (in general) of strictly positive dimension, and furthermore are (in general) not proper if $\calG\to X$ is not parahoric.

\end{document}